\documentclass[12pt]{article}
\usepackage{graphics, graphicx}     %%%   Ïàêåòû äëÿ ðàáîòû ñ ðèñóíêàìè.

\usepackage{cite}
\usepackage{amsmath}
\usepackage{amssymb}
\usepackage{amsthm}

\theoremstyle{plain}
\newtheorem{Theorem}{Theorem}[section]  %
\newtheorem{Lemma}{Lemma}[section] 
\newtheorem{Proposition}{Proposition}[section] 
\theoremstyle{definition}
\newtheorem{Remark}{Remark}[section]

\theoremstyle{definition}
\newenvironment{Proof} % èìÿ îêðóæåíèÿ
{\par\noindent{\it Proof of}} % êîìàíäû äëÿ \begin
{\hfill$\vspace{5mm}\scriptstyle\blacksquare$} % êîìàíäû äëÿ \end

\numberwithin{equation}{section} %×òîáû íóìåðàöèÿ â êàæäîé ñåêöèè áûëà íåçàâèñèìîé
\numberwithin{figure}{section} %×òîáû íóìåðàöèÿ â êàæäîé ñåêöèè áûëà íåçàâèñèìîé
\numberwithin{table}{section} %×òîáû íóìåðàöèÿ â êàæäîé ñåêöèè áûëà íåçàâèñèìîé

\begin{document}

\setcounter{page}{1}

\markboth{M.I. Isaev}{Asymptotic enumeration of Eulerian circuits for graphs with strong mixing properties}

\title{Asymptotic enumeration of Eulerian circuits for graphs with strong mixing properties}
\date{}
\author{ {\bf M.I. Isaev } }
     
\maketitle
%{\bf Abstract}
\begin{abstract}
 	We prove an asymptotic formula for the number of Eulerian circuits for graphs with strong
		mixing properties and with vertices having even degrees. The exact value is determined up to the multiplicative error $O(n^{-1/2+\varepsilon})$, where $n$ is the number of vertices.
\end{abstract}

\section{Introduction}
$\ \ \ $ Let $G$ be a simple connected graph all of whose vertices have even degrees. An {\it Eulerian
circuit } in $G$ is a closed walk (see, for example, \cite{Biggs1976}) which uses every edge of $G$ exactly once. 
Two Eulerian circuits are called equivalent if one is a cyclic permutation of the other. 
It is clear that the size of such an equivalence class equals the
number of edges of graph $G$. Let $EC(G)$ denote the number of equivalence classes
of Eulerian circuits in $G$.

	The problem of counting the number of  Eulerian circuits in an undirected simple graph (graph without loops and multiple edges) is
complete for the class $\# P$, i.e. the existence of a polynomial algorithm for this problem implies the existence of
a polynomial algorithm for any problem in the class $\# P$ and, in particular, the equivalence of the classes $P$ and $NP$. 
 (see \cite{Brightwell2005}). In other words, the problem of counting the number of  Eulerian circuits is difficult in terms of the complexity theory. Moreover, it should be noted that in contrast to many other hard problems of counting on graphs (see, for example, \cite{AKK1995}, \cite{Mihail1996}), even approximate and probabilistic polynomial algorithms for counting the number of Eulerian circuits have not been obtained for the general case and are known only for some special classes of graphs with low density, see \cite{CCM2012} and \cite{TV2001}.

	As concerns the class of complete graphs $K_n$,  the exact expression of the number of Eulerian circuits for odd $n$ is unknown (it is clear that $EC(K_n)=0$ for even $n$) and 	only the asymptotic formula was obtained (see \cite{Brendan1995}): 
	as $n\rightarrow \infty$ with $n$ odd
	\begin{equation}\label{Eq_1_1}
		\begin{aligned}
		EC(K_n) =  2^{\frac{(n-1)^2}{2}} \pi^{-\frac{n-1}{2}} n^{\frac{n-2}{2}}
		  \left(\left( \frac{n-1}{2} -1 \right)!\right)^n \Big(1 + O(n^{-1/2+\varepsilon})\Big)=\\
		  = 2^{\frac{n+1}{2}} \pi^{\frac{1}{2}}\, e^{-\frac{n^2}{2}+\frac{11}{12}}\, n^{\frac{(n-2)(n+1)}{2}}
		  \Big(1 + O(n^{-1/2+\varepsilon})\Big)
		\end{aligned}
	\end{equation}
	for any fixed $\varepsilon>0$.
	
In \cite{Isaev2011} the analytic approach 
of \cite{Brendan1995} was generalized. This approach is based on expression of the result in terms of a multidimensional integral and its estimation as dimension tends to infinity. 
 In particular, the asymptotic 
behaviour of the number of Eulerian circuits was determined for graphs with large algebraic connectivity. 
This class of graphs we mean as the class of graphs having strong mixing properties.

There are several classic graph parameters which express mixing properties of graphs:
  the algebraic connectivity, the Cheeger constant (isoperimetric number), the spectral gap between the $1$ and the second largest eigenvalue of the transition probability of the random walk on a graph.
  It should be noted that, using any of these parametres, one can get equivalent definitions of the class of graphs having strong mixing properties (for more detailed information, see \cite{IK2012}).
 
In addition, it is shown in \cite{IK2012} that a random graph in the Gilbert graph model $G(n,p)$ 
(each possible edge occurs in a graph with $n$ vertices independently with probability $p$) for  $n\rightarrow \infty$
and fixed $p>0$ has strong mixing properties with probability close to $1$ (with the exception of some exponentially small of $n$ value).

In the present work we continue studies of \cite{Isaev2011}, \cite{IK2012}, \cite{Brendan1995}. We prove an asymptotic formula for the
number of Eulerian circuits of graphs having strong mixing properties. This result 
 is presented in detail in Section 2 of the present work.

Actually, the estimation of the number of Eulerian circuits
was reduced in  \cite{Isaev2011} to estimating of an $n$-dimensional integral which is close to Gaussian-type. 
We partly  repeat this reduction  in Sections 3, 8 of the present paper. In addition, 
 we develop an approach for estimating of integrals of such a type in Sections 4, 6, 7.  We prove the main result in Section 5.

An orientation of edges of a graph 
such that at each vertex the number of incoming and outgoing edges are the same is  called {\ it
Eulerian orientation}.
	In \cite{Isaev2012} the asymptotic 
behaviour of the number of Eulerian orientations was determined for graphs having strong mixing properties. 
Apparently, proceeding from the results of \cite{Isaev2012} and the estimates of the present work, it is possible to prove the asymptotic formula for the number of Eulerian orientations given in \cite{IK2012}. In a subsequent paper we plan to develop this approach. 
%%%%%%%%%%%%%%%%%%%%%%%%%%%%%%%%%%%%%%%%%%%%%%%%%%%%%%%%%%%%%%%%%%%%%%%%%%%%%%%%%%%%%%%%%%%%%%%%%%%%%%%%%%%%%%%%%%%%%%%%%
%%%%%%%%%%%%%%%%%%%%%%%%%%%%%%%%%%%%%%%%%%%%%%%%%%%%%%%%%%%%%%%%%%%%%%%%%%%%%%%%%%%%%%%%%%%%%%%%%%%%%%%%%%%%%%%%%%%%%%%%%
%%%%%%%%%%%%%%%%%%%%%%%%%%%%%%%%%%%%%%%%%%%%%%%%%%%%%%%%%%%%%%%%%%%%%%%%%%%%%%%%%%%%%%%%%%%%%%%%%%%%%%%%%%%%%%%%%%%%%%%%%
%%%%%%%%%%%%%%%%%%%%%%%%%%%%%%%%%%%%%%%%%%%%%%%%%%%%%%%%%%%%%%%%%%%%%%%%%%%%%%%%%%%%%%%%%%%%%%%%%%%%%%%%%%%%%%%%%%%%%%%%%
%%%%%%%%%%%%%%%%%%%%%%%%%%%%%%%%%%%%%%%%%%%%%%%%%%%%%%%%%%%%%%%%%%%%%%%%%%%%%%%%%%%%%%%%%%%%%%%%%%%%%%%%%%%%%%%%%%%%%%%%%
%%%%%%%%%%%%%%%%%%%%%%%%%%%%%%%%%%%%%%%%%%%%%%%%%%%%%%%%%%%%%%%%%%%%%%%%%%%%%%%%%%%%%%%%%%%%%%%%%%%%%%%%%%%%%%%%%%%%%%%%%
%%%%%%%%%%%%%%%%%%%%%%%%%%%%%%%%%%%%%%%%%%%%%%%%%%%%%%%%%%%%%%%%%%%%%%%%%%%%%%%%%%%%%%%%%%%%%%%%%%%%%%%%%%%%%%%%%%%%%%%%%
%%%%%%%%%%%%%%%%%%%%%%%%%%%%%%%%%%%%%%%%%%%%%%%%%%%%%%%%%%%%%%%%%%%%%%%%%%%%%%%%%%%%%%%%%%%%%%%%%%%%%%%%%%%%%%%%%%%%%%%%%
\section{Main result}
Let $G$ be an undirected simple graph with vertex set $VG=\{v_1,v_2,\ldots,v_n\}$ and edge set $EG$.
	We define $n\times n$ matrix $Q$ by
	\begin{equation}\label{Laplace_def}
		Q_{jk} = 
		\left\{
			\begin{array}{cl}
			-1, & \{v_j,v_k\}\in EG,\\
			\phantom{-}d_j,& j = k,\\
			\phantom{-}0,  & \text{ otherwise,}
		\end{array}\right.
\end{equation}
	where $n = |VG|$ and  $d_j$ denotes the degree of $v_j\in VG$. The matrix $Q = Q(G)$ is called the {\it Laplacian matrix} 
	of the graph $G$. The eigenvalues 
	$\lambda_1 \leq \lambda_2 \leq \ldots \leq \lambda_n$ of the
matrix $Q$ are always non-negative real numbers and the number of zero eigenvalues of $Q$ coincides
with the number of connected components of $G$, in particular, $\lambda_1 = 0$. The eigenvalue $\lambda_2 = \lambda_2(G)$ is called the
{\it algebraic connectivity} of the graph $G$. In addition, the following inequalities hold:
\begin{equation}\label{d_lambda}
	2\min_j d_j - n + 2 \leq \lambda_2 \leq \frac{n}{n-1} \min_j d_j.
\end{equation}
For more information on the spectral properties
of graphs see, for example, \cite{Fiedler1973} and \cite{Mohar1991}. 

An 	acyclic connected subgraph of the graph $G$d which includes all of its vertices is called a {\it spanning tree} of $G$.
	According to Kirchhoff's Matrix-Tree-Theorem, see \cite{Kirchoff1847}, we have that
\begin{equation}\label{Eq_1_3}
	t(G) = \frac{1}{n}\lambda_2\lambda_3\cdots\lambda_{n} = \det{M_{11}},
\end{equation}
where $t(G)$ denotes the number of spanning trees of the graph $G$ and $M_{11}$
results from deleting the first row and the first column of $Q$.
	
	We call the graph $G$ as $\gamma$-mixing graph, $\gamma>0$,  if 
	\begin{equation}\label{ac_cond}
		\text{the algebraic connectivity } \lambda_2 = \lambda_2(G) \geq \gamma |VG|.
	\end{equation}
The main result of the present work is the following theorem.
\begin{Theorem}\label{main}
	Let $G$ be an undirected simple graph  with  $n$ vertices $v_1, v_2, \ldots, v_n$ having even degrees. Let $G$ be a $\gamma$-mixing graph 
	for some $\gamma>0$. 
	%\begin{equation}\label{condition_main}
	%	\left\|(Q+J)^{-1}\right\|_{\infty} \leq \frac{1}{\sigma n},
	%\end{equation}  
	%where $J$ denotes the matrix with every entry $1$. 
%	for some $\sigma > 1/2$ the degree of each vertex of $G$ at least $\sigma n$.
	Then 
	\begin{equation}\label{main_eq}
		\begin{aligned}
		EC(G) = \left(1 + \delta(G) \right) e^{K_{ec}} 
		\left(
		   		   2^{|EG|-\frac{n-1}{2}} \pi^{-\frac{n-1}{2}} { \sqrt{t(G)}}\,
		   		   \prod\limits_{j=1}\limits^{n}\left(\frac{d_j}{2}-1\right)! 
		\right),\\
		K_{ec} = - \frac{1}{4} \sum\limits_{\{v_j,v_k\}\in EG} \left(\frac{1}{d_j+1} - \frac{1}{d_k+1}\right)^2,
			\end{aligned}
	\end{equation}
	where  $EG$ denotes the  edge set of $G$, $d_j$ is  the degree of vertex $v_j$,	  $t(G)$ is the number of spanning trees of  $G$ 
and for any  $\varepsilon>0$ 
	\begin{equation}\label{main_eq2}
		|\delta(G)| \leq C n^{-1/2+\varepsilon},
	\end{equation}
	where constant $C > 0$ depends only on $\gamma$ and $\varepsilon$.   
\end{Theorem}

Proof of Theorem \ref{main} is given in Section 5. This proof is based on results presented in
Sections 3, 4.

%\noindent
%{\bf Remark 2.1.} 

\begin{Remark}
For the case of the complete graph we have that: 
\begin{equation}\label{Remark_eq1}
\lambda_2(K_n) = n,\  EK_n = \frac{n(n-1)}{2},\  t(K_n) = n^{n-2},\  
K_{ec} = 0.
\end{equation} 
We obtain that the result of Theorem \ref{main} for the case of the complete graph is	equivalent to (\ref{Eq_1_1}).
\end{Remark}	
%%%%%%%%%%%%%%%%%%%%%%%%%%%%%%%%%%%%%%%%%%%%%%%%%%%%%%%%%%%%%%%%%%%%%%%%%%%%%%%%%%%%%%%%%%%%%%%%%%%%%%%%%%%%%%%%%%%%%%%%%
%%%%%%%%%%%%%%%%%%%%%%%%%%%%%%%%%%%%%%%%%%%%%%%%%%%%%%%%%%%%%%%%%%%%%%%%%%%%%%%%%%%%%%%%%%%%%%%%%%%%%%%%%%%%%%%%%%%%%%%%%
%%%%%%%%%%%%%%%%%%%%%%%%%%%%%%%%%%%%%%%%%%%%%%%%%%%%%%%%%%%%%%%%%%%%%%%%%%%%%%%%%%%%%%%%%%%%%%%%%%%%%%%%%%%%%%%%%%%%%%%%%
%%%%%%%%%%%%%%%%%%%%%%%%%%%%%%%%%%%%%%%%%%%%%%%%%%%%%%%%%%%%%%%%%%%%%%%%%%%%%%%%%%%%%%%%%%%%%%%%%%%%%%%%%%%%%%%%%%%%%%%%%
%%%%%%%%%%%%%%%%%%%%%%%%%%%%%%%%%%%%%%%%%%%%%%%%%%%%%%%%%%%%%%%%%%%%%%%%%%%%%%%%%%%%%%%%%%%%%%%%%%%%%%%%%%%%%%%%%%%%%%%%%
%%%%%%%%%%%%%%%%%%%%%%%%%%%%%%%%%%%%%%%%%%%%%%%%%%%%%%%%%%%%%%%%%%%%%%%%%%%%%%%%%%%%%%%%%%%%%%%%%%%%%%%%%%%%%%%%%%%%%%%%%
%%%%%%%%%%%%%%%%%%%%%%%%%%%%%%%%%%%%%%%%%%%%%%%%%%%%%%%%%%%%%%%%%%%%%%%%%%%%%%%%%%%%%%%%%%%%%%%%%%%%%%%%%%%%%%%%%%%%%%%%%
\section{Reduction to the integral}

A {\it directed tree with root } $v$ is a connected directed graph $T$ such that $v \in VT$ has out-
degree zero, and each other vertex has out-degree one. Thus, $T$ is a tree which has each edge
oriented towards $v$.

	Let $G$ be a connected undirected simple graph  with  $n$ vertices $v_1, v_2, \ldots, v_n$ having even degrees.
	Note that for every spanning tree $T$ of the graph $G$ and any vertex $v_r\in VG$ there is only one orientation of the edges of $T$
 such that we obtain a directed tree with root $v_r$. We denote by ${\cal T}_r$ the set of directed trees with root $v_r$ obtained in such a way.
 
 We recall that (see Section 4 and formulas (4.6), (4.7) of \cite{Isaev2011}):
\begin{equation}\label{EulS}
EC(G) = \prod\limits_{j=1}^n \big( \frac{d_j}{2} -1 \big)! \, 2^{|EG|-n+1} \pi^{-n} S,
\end{equation}
where for any $r\in \mathbb{N}$,  $r\leq n$,
\begin{equation}\label{S}
	S = \int\limits_{U_n(\pi/2)} 
 % \left( 
     \prod\limits_{\{v_j,v_k\}\in EG} \cos \Delta_{jk} \sum \limits_{T \in {\cal T}_r} 
     \prod \limits_{(v_j,v_k)\in ET} (1+i\tan \Delta_{jk})  
 % \right) 
 \
 d\vec{\xi},
\end{equation} 
where $\Delta_{jk} = \xi_j-\xi_k$ and
\begin{equation} 
U_n(\rho) = \{(\xi_1, \xi_2, \ldots, \xi_n)\in \mathbb{R}^n\ : \ |\xi_j| \leq \rho \text{ for all } j= {1\ldots n}\}.
\end{equation}

	We approach the integral by first estimating it in the region which is  the asymptotically significant one. In what follows, 
	we fix some small constant $\varepsilon>0$. Define 
\begin{equation}
	\begin{aligned}
	V_0 = \{\vec{\xi}\in U_n(\pi/2) :\ |\xi_j - \xi_k|_\pi \leq n^{-1/2+\varepsilon} \text { for any } 
	1\leq j,k\leq n \}\\
	|\xi_j - \xi_k|_\pi = \min\limits_{l\in \mathbb{Z}} |\xi_j - \xi_k+ \pi l|.
	\end{aligned}
\end{equation}and let $S_0$ denote the contribution  to $S$ of $\vec{\xi}\in V_0$:
\begin{equation}
	S_0 =  \frac{1}{n} \sum \limits_{r=1}\limits^n \int\limits_{U_n(\pi/2)} 
 % \left( 
     \prod\limits_{\{v_j,v_k\}\in EG} \cos \Delta_{jk} \sum \limits_{T \in {\cal T}_r} 
     \prod \limits_{(v_j,v_k)\in ET} (1+i\tan \Delta_{jk})  
 % \right) 
 \
 d\vec{\xi}.
\end{equation} 

In this section  we use standart notation $f = O(g)$ as $n\rightarrow \infty$ 
which indicates that there exist $c,n_0>0$ such that for $n\geq n_0$ the inequality $|f| \leq c|g|$ holds.. 

Under assumptions of Theorem \ref{main}, we have that as $n\rightarrow \infty$ 
	\begin{equation}\label{S+S_0}
		S  = \left(1 + O\left(\exp(-cn^{2\varepsilon})\right)\right) S_0
	\end{equation}
	for some $c>0$ depending only on $\gamma$. For the proof of \eqref{S+S_0}, see Theorem 6.3 of \cite{Isaev2011}.

Let
\begin{equation}\label{hatQ}
	W = \hat{Q}^{-1} = (Q + J)^{-1},
\end{equation}
where $Q$ is the Laplacian matrix and $J$ denotes the matrix with every entry $1$.
Let $\vec{\alpha} = (\alpha_1,\ldots,\alpha_n)\in \mathbb{R}^n$ be defined by
\begin{equation}\label{alpha}
 	\alpha_j = W_{jj}.
\end{equation}
Let 
\begin{equation}\label{R_theta}
 R(\vec{\xi}) = \mbox{tr} (\Lambda(\vec{\xi}) W \Lambda(\vec{\xi}) W),
\end{equation}
where $\mbox{tr}(\cdot)$ is the trace fucntion, $\Lambda(\vec{\xi})$ denotes the diagonal
matrix whose diagonal elements are equal to corresponding  components of the vector $Q\vec{\xi}$.

The sum over ${\cal T}_r$ in the integrand of (\ref{S}) can be expressed as a determinant, according to the following theorem of \cite{Tutte1948}, which is a generalization of aforementioned  Kirchhoff's Matrix-Tree-Theorem:

\begin{Theorem}\label{Tutte}
Let $w_{jk}$ $(1 \leq j, k \leq n$, $j \neq k$) be arbitrary. Define the $n \times n$ matrix $A$ by
\begin{equation}\label{matrix_Tutte}
A_{jk} = 
\left\{
\begin{array}{cc}
-w_{jk},& \text{ if } j\neq k,\\
\sum_{r\neq j} w_{jr},& \text{ if } k=j
\end{array}
\right.,
\end{equation}
the sum being over $1 \leq r \leq n$ with $r\neq j$. For any $r$ with $1 \leq r \leq n$, let 
 $M_r$ denote the
principal minor of $A$ formed by removing row $r$ and column $r$. Then
\begin{equation}
\det M_r = \sum\limits_{T} \prod\limits_{(v_j, v_k)\in ET} w_{jk},
\end{equation}
where the sum is over all directed trees $T$ with $VT = \{v_1, v_2, \ldots, v_n\}$ and root $v_r$.
\end{Theorem}

Using formulas \eqref{S} and \eqref{S+S_0}, Theorem \ref{Tutte} and the Taylor series expansion of $\cos \Delta_{jk}$ and 
$\tan \Delta_{jk}$in the region $V_0$, one can obtain the following proposition:

\begin{Proposition}\label{Proposition_3.1}
	Let the assumptions of Theorem \ref{main} hold. Then as $n\rightarrow \infty$
	\begin{equation}\label{S_0+++}
		S_0  = \left(1 + O\left(n^{-1/2+6\epsilon}\right)\right) 2^{-1/2}\pi^{1/2} n^{-1}  \det{\hat{Q}} \ {\rm Int},
	\end{equation}
	\begin{equation}
		{\rm Int} =\int\limits_{ U_n(n^{-1/2+\varepsilon})}	
		\exp\Bigg(
				i\,\vec{\xi}^T Q \vec{\alpha} - \frac{1}{2} \vec{\xi}^T \hat{Q} \vec{\xi} 
			- \frac{1}{12}\sum\limits_{\{v_j,v_k\}\in EG} \Delta_{jk}^4 + \frac{1}{2}R(\vec{\xi})
			\Bigg) 
		 d\vec{\xi},
	\end{equation}
	where $\hat{Q}$, $\vec{\alpha}$ and $R(\vec{\xi})$ are the same as in \eqref{hatQ}, \eqref{alpha} and \eqref{R_theta}, respectively.
\end{Proposition}

We prove in detail  Proposition \ref{Proposition_3.1} in Section 8. Actually, this proof was implicitly given in \cite{Isaev2011} 
(see Lemma 5.3 of \cite{Isaev2011}).

Thus, we get that to prove Theorem \ref{main} it remains only to estimate the integral $\mbox{Int}$ of \eqref{S_0+++}.
%%%%%%%%%%%%%%%%%%%%%%%%%%%%%%%%%%%%%%%%%%%%%%%%%%%%%%%%%%%%%%%%%%%%%%%%%%%%%%%%%%%%%%%%%%%%%%%%%%%%%%%%%%%%%%%%%%%%%%%%%
%%%%%%%%%%%%%%%%%%%%%%%%%%%%%%%%%%%%%%%%%%%%%%%%%%%%%%%%%%%%%%%%%%%%%%%%%%%%%%%%%%%%%%%%%%%%%%%%%%%%%%%%%%%%%%%%%%%%%%%%%
%%%%%%%%%%%%%%%%%%%%%%%%%%%%%%%%%%%%%%%%%%%%%%%%%%%%%%%%%%%%%%%%%%%%%%%%%%%%%%%%%%%%%%%%%%%%%%%%%%%%%%%%%%%%%%%%%%%%%%%%%
%%%%%%%%%%%%%%%%%%%%%%%%%%%%%%%%%%%%%%%%%%%%%%%%%%%%%%%%%%%%%%%%%%%%%%%%%%%%%%%%%%%%%%%%%%%%%%%%%%%%%%%%%%%%%%%%%%%%%%%%%
%%%%%%%%%%%%%%%%%%%%%%%%%%%%%%%%%%%%%%%%%%%%%%%%%%%%%%%%%%%%%%%%%%%%%%%%%%%%%%%%%%%%%%%%%%%%%%%%%%%%%%%%%%%%%%%%%%%%%%%%%
%%%%%%%%%%%%%%%%%%%%%%%%%%%%%%%%%%%%%%%%%%%%%%%%%%%%%%%%%%%%%%%%%%%%%%%%%%%%%%%%%%%%%%%%%%%%%%%%%%%%%%%%%%%%%%%%%%%%%%%%%
%%%%%%%%%%%%%%%%%%%%%%%%%%%%%%%%%%%%%%%%%%%%%%%%%%%%%%%%%%%%%%%%%%%%%%%%%%%%%%%%%%%%%%%%%%%%%%%%%%%%%%%%%%%%%%%%%%%%%%%%%
%%%%%%%%%%%%%%%%%%%%%%%%%%%%%%%%%%%%%%%%%%%%%%%%%%%%%%%%%%%%%%%%%%%%%%%%%%%%%%%%%%%%%%%%%%%%%%%%%%%%%%%%%%%%%%%%%%%%%%%%%
%%%%%%%%%%%%%%%%%%%%%%%%%%%%%%%%%%%%%%%%%%%%%%%%%%%%%%%%%%%%%%%%%%%%%%%%%%%%%%%%%%%%%%%%%%%%%%%%%%%%%%%%%%%%%%%%%%%%%%%%%
\section{Asymptotic estimates of integrals}
We fix constants $a,b,\varepsilon>0$. In this section we use notation $f = O(g)$ meaning that $|f|\leq c|g|$ for some 
$c>0$ depending only on $a,b$ and $\varepsilon$. 

%\subsection{ Conventions and notations}

Let $p \geq 1$ be a real number and $\vec{x}\in \mathbb{R}^n$. Let
\begin{equation}
	\left\|\vec{x}\right\|_p = \left(\sum\limits_{j=1}^{n} |x_j|^p \right)^{1/p}.
\end{equation}
For $p = \infty$ we have the maximum norm 
\begin{equation}
	\left\|\vec{x}\right\|_\infty = \max_j|x_j|.
\end{equation}
The matrix norm corresponding to the $p$-norm for vectors is
\begin{equation}
	\left\|A\right\|_p = \sup_{\vec{x}\neq 0} \frac{\left\|A\vec{x}\right\|_p}{\left\|\vec{x}\right\|_p}.
\end{equation}
One can show that for symmetric matrix $A$ and $p\geq 1$
\begin{equation}\label{Spec_norm}
	\left\|A\right\|_p \geq \left\|A\right\|_2.
\end{equation}
For invertible matrices define the condition number
\begin{equation}
	\mu_p(A) = \left\|A\right\|_p\cdot\left\|A^{-1}\right\|_p \geq \left\|A A^{-1}\right\|_p = 1.
\end{equation}
%Îáîçíà÷èì $\left\|A\right\|_{HS}$ íîðìó Ãèëüáåðòà-Øìèäòà ìàòðèöû $A$.
%\begin{equation}
% \left\|A\right\|_{HS} = \sqrt{\sum\limits_{j=1}\limits^{n}\sum\limits_{k=1}\limits^{n} |A_{jk}|^2}
%\end{equation}

Let $I$ be identity $n\times n$ matrix and $A=I+X$ be such a matrix that:
\begin{equation}\label{A_ass}
	\begin{aligned}
		\text{ $A$ is positive definite symmetric matrix},\\
		|X_{jk}| \leq a/n, \ \ X_{jj} = 0, \ \ 
		\|A^{-1}\|_2 \leq b.
	\end{aligned}
\end{equation}
Note that 
\begin{equation}\label{3.7}
\|A^{-1}\|_2^{-1} \leq \|A\|_2 \leq \|{A}\|_\infty = \|{A}\|_1 = \max_{j}{\sum\limits_{k=1}^{n}} |{A}_{jk}| = O(1).
\end{equation}
We recall that (see Lemma 3.2 of \cite{Isaev2012} ), under assumptions (\ref{A_ass}), 
\begin{equation}\label{3.8}
	\mu_\infty(A) = \mu_1(A) = O(\mu_2(A)).
\end{equation} 
Using (\ref{Spec_norm}), (\ref{A_ass}), (\ref{3.7}) and (\ref{3.8}), we obtain the following lemma:
\begin{Lemma}\label{Lemma_A}
	Let $A$ satisfy (\ref{A_ass}). Then 
	\begin{equation}\label{3.9}
			\|{A}^{-1}\|_\infty = \|{A}^{-1}\|_1 = O(1),
	\end{equation}
	\begin{equation}\label{3.10}
			|X'_{jk}| = O(n^{-1}),
	\end{equation}
	where
	\begin{equation}
		X' = A^{-1} - I = A^{-1}(I - A) = - A^{-1} X.
	\end{equation}
\end{Lemma}
%\begin{Proof} { \it Lemma \ref{Lemma_A}.}
%	sd
%\end{Proof}
%\subsection{An integral}

We use the following notation:
\begin{equation}
		<g>_{F,\Omega} = \int\limits_{\Omega} g(\vec{\theta}) e^{F(\vec{\theta})} d\vec{\theta}, 	
\end{equation}
where $g,F$ are some functions on $\mathbb{R}^n$.
For $r>0$ let
\begin{equation}
	<g>_{F,r} = <g>_{F,U_n(rn^{\varepsilon})},
\end{equation}
where 
\begin{equation} 
U_n(\rho) = \{(\theta_1, \theta_2, \ldots, \theta_n)\in \mathbb{R}^n\ : \ |\theta_j| \leq \rho \text{ for all } j= {1\ldots n}\}.
\end{equation} 
%In the case of $r=1$ we denote $<g>_{F} = <g>_{F,1}$. 
We use functions  of the following type:
\begin{equation}\label{def_F}
	F(\vec{\theta}) = - \vec{\theta}^T A \vec{\theta} + H(\vec{\theta}),
\end{equation}
where $A$ satisfy (\ref{A_ass}). Let consider the following assumptions on function $H$, which we will need further:
	\begin{equation}\label{R_1}	
		H(\vec{\theta}) \leq c_1 \frac{\vec{\theta}^T A \vec{\theta}}{n},
	\end{equation}
	\begin{equation}\label{R_2}
		\displaystyle
		 \left\|\frac{\partial H(\vec{\theta})}{\partial \vec{\theta}}\right\|_\infty \leq c_2 
		 \frac{\|\vec{\theta}\|_\infty^3 + \|\vec{\theta}\|_\infty}{n} 
	 \end{equation}
%	for any $\vec{\theta} \in U_n(rn^{\varepsilon})$. 
	%where constants $c_1,c_2>0$ depend only on $r$, $a$, $b$ and $\varepsilon$. 	
	For the case when $H \equiv 0$ we use notations:  
\begin{equation}
<g>_\Omega =<g>_{F,\Omega},\  \ <g>_r =<g>_{F,r},\ \ 
<g> =<g>_{+\infty}.
\end{equation}
%Let $\vec{\phi}(\vec{\theta}) = A \vec{\theta}$. 

%\begin{Proposition}
%	Let $A$ satisfy (\ref{A_ass}) and  assumption (\ref{R_1}) holds. Let 
%	\begin{equation}
%	|T(\vec{\theta})| \leq c (\|\vec{\theta}\|_\infty+1)^m
%	\end{equation}
%	for some fixed $c,m>0$. Then
%		\begin{equation}
%		<T(\vec{\theta})>_{F,r_1} - <T(\vec{\theta)}>_{F,r_2}  =  O\left(\exp(-c_3n^{2\varepsilon})\right)<1>,
%	\end{equation} 
%	where $F$ is defined by (\ref{def_F}) and constant $c_3>0$ depends only on $c$, $m$, $r_1$, $r_2$, $a$, $b$ and $\varepsilon$.
%\end{Proposition}
\begin{Lemma}\label{Proposition_4.1}
Let $\Omega\subset \mathbb{R}^n$ be such that  $U_n(r_1n^{\varepsilon}) \subset \Omega \subset U_n(r_2n^{\varepsilon})$ for some $r_2>r_1>0$. 
Let $A$ satisfy (\ref{A_ass}) and  assumptions (\ref{R_1}), (\ref{R_2}) hold for some $c_1,c_2>0$. Then
	\begin{equation}\label{Omega_infty}
		<1>_\Omega = (1 + O\left(\exp(-c_3n^{2\varepsilon})\right))<1>,
	\end{equation}
	\begin{equation}\label{theta^0}
		<1>_{F,\Omega} = O\left(<1>\right),
	\end{equation}
	\begin{equation}\label{theta^2}
		<\theta_k^2>_{F,\Omega} = \frac{1}{2}<1>_{F,\Omega} + O(n^{-1+4\varepsilon})<1>,
	\end{equation}
	\begin{equation}\label{theta^4}
		<\theta_k^4>_{F,\Omega} = \frac{3}{4}<1>_{F,\Omega} + O(n^{-1+7\varepsilon})<1>
	\end{equation}
	and,  for $k\neq l$:
	\begin{equation}\label{theta^11}
		<\theta_k\theta_l>_{F,\Omega} = O(n^{-1+5\varepsilon})<1>,
	\end{equation}
	\begin{equation}\label{theta^13}
		<\theta_k\theta_l^3>_{F,\Omega} = O(n^{-1+7\varepsilon})<1>,
	\end{equation}
	\begin{equation}\label{theta^22}
		<\theta_k^2\theta_l^2>_{F,\Omega} = \frac{1}{4}<1>_{F,\Omega} + O(n^{-1+7\varepsilon})<1>,
	\end{equation}
	where $F$ is defined by (\ref{def_F}) and $c_3 = c_3(r_1, r_2, c_1,c_2, a, b, \varepsilon)>0$.
	
	In addition, for any vector $\vec{p} = (p_1, p_2, \ldots p_n) \in \mathbb{R}^n$, $\|\vec{p}\|_\infty = O(n^{-1/2})$,    
	\begin{equation}\label{theta+o}
		\begin{aligned}
		<\theta_k e^{i \vec{\theta}^T\vec{p} - i p_k\theta_k}>_{F,\Omega} &= \\=
	 		\frac{i}{2} \sum\limits_{j\neq k, j\leq n} p_j (A^{-1})_{jk} 
	 		< &e^{i\vec{\theta}^T\vec{p} - i p_k\theta_k}>_{(F - \frac{1}{2}p_k^2 \theta_k^2), \Omega}+\\ &+ O(n^{-1+5\varepsilon})<1>,
	 		\end{aligned}
	\end{equation}
	where  $(A^{-1})_{jk}$ denotes $(j,k)$-th element of the matrix $A^{-1}$.
\end{Lemma}
Proof of Lemma \ref{Proposition_4.1} is given in Section 6.

%%%%%%%%%%%%%%%%%%%%%%%%%%%%%%%%%%%%%%%%%%%%%%%%%%%%%%%%%%%%%%%%%%%%%%%%%%%%%%%%%%%%%%%%%%%%%%%%%%%%%%%%%%%%%%%%%%%%%%%%%
%%%%%%%%%%%%%%%%%%%%%%%%%%%%%%%%%%%%%%%%%%%%%%%%%%%%%%%%%%%%%%%%%%%%%%%%%%%%%%%%%%%%%%%%%%%%%%%%%%%%%%%%%%%%%%%%%%%%%%%%%
%%%%%%%%%%%%%%%%%%%%%%%%%%%%%%%%%%%%%%%%%%%%%%%%%%%%%%%%%%%%%%%%%%%%%%%%%%%%%%%%%%%%%%%%%%%%%%%%%%%%%%%%%%%%%%%%%%%%%%%%%
%%%%%%%%%%%%%%%%%%%%%%%%%%%%%%%%%%%%%%%%%%%%%%%%%%%%%%%%%%%%%%%%%%%%%%%%%%%%%%%%%%%%%%%%%%%%%%%%%%%%%%%%%%%%%%%%%%%%%%%%%
%%%%%%%%%%%%%%%%%%%%%%%%%%%%%%%%%%%%%%%%%%%%%%%%%%%%%%%%%%%%%%%%%%%%%%%%%%%%%%%%%%%%%%%%%%%%%%%%%%%%%%%%%%%%%%%%%%%%%%%%%
%%%%%%%%%%%%%%%%%%%%%%%%%%%%%%%%%%%%%%%%%%%%%%%%%%%%%%%%%%%%%%%%%%%%%%%%%%%%%%%%%%%%%%%%%%%%%%%%%%%%%%%%%%%%%%%%%%%%%%%%%
%%%%%%%%%%%%%%%%%%%%%%%%%%%%%%%%%%%%%%%%%%%%%%%%%%%%%%%%%%%%%%%%%%%%%%%%%%%%%%%%%%%%%%%%%%%%%%%%%%%%%%%%%%%%%%%%%%%%%%%%%
%%%%%%%%%%%%%%%%%%%%%%%%%%%%%%%%%%%%%%%%%%%%%%%%%%%%%%%%%%%%%%%%%%%%%%%%%%%%%%%%%%%%%%%%%%%%%%%%%%%%%%%%%%%%%%%%%%%%%%%%%
%%%%%%%%%%%%%%%%%%%%%%%%%%%%%%%%%%%%%%%%%%%%%%%%%%%%%%%%%%%%%%%%%%%%%%%%%%%%%%%%%%%%%%%%%%%%%%%%%%%%%%%%%%%%%%%%%%%%%%%%%
\section{Proof of Theorem \ref{main}}

The Laplacian matrix $Q$  of the graph $G$ ( 
defined in \eqref{Laplace_def}) has the
eigenvector $[1,1,\ldots,1]^T$, corresponding to the eigenvalue $\lambda_1 = 0$.
%For our purpose it is convenient to denote by $L$ the orthogonal complement to $[1,1,\ldots,1]^T$. 
Let $\hat{Q} = Q + J$, where  $J$ denotes the matrix with every entry $1$. Note that $Q$ and $\hat{Q}$ 
have the same set of
eigenvectors and eigenvalues, except for the eigenvalue corresponding to the eigenvector $[1,1,\ldots,1]^T$, which equals 
$0$ for $Q$ and $n$ for $\hat{Q}$. Using (\ref{Eq_1_3}), we find that
\begin{equation}\label{tGQ}
	t(G) = \frac{1}{n}\lambda_2\lambda_3\cdots\lambda_{n} = \frac{\det{\hat{Q}}}{n^2}.
\end{equation}
Using (\ref{Spec_norm}), we get that
\begin{equation}\label{eq5.3}
	\lambda_{n} = ||Q||_2 \leq ||\hat{Q}||_2 \leq ||\hat{Q}||_1 = \max_{j}{\sum\limits_{k=1}^{n}} |\hat{Q}_{jk}| = n.
\end{equation}
Then, we find that
\begin{equation}\label{4.2}
	||\hat{Q}^{-1}||_2 = \frac{1}{\lambda_2} \leq \frac{1}{\gamma n}.
\end{equation}
Using (\ref{d_lambda}), we get that
\begin{equation}\label{4.3} 
	  n-1 \geq d_j \geq \lambda_2 \frac{n-1}{n}\geq \gamma (n-1),
\end{equation}
where $d_j$ is the degree of  $v_j$.
Consider  the integral of Proposition \ref{Proposition_3.1}:
\begin{equation}\label{eq5.5}
	\begin{aligned}
		\mbox{Int} =\int\limits_{ U_n(n^{-1/2+\varepsilon})}	
		\exp\Bigg(
				i\,\vec{\xi}^T \vec{\beta} - \frac{1}{2} \vec{\xi}^T \hat{Q} \vec{\xi} 
			- \frac{1}{12}\sum\limits_{\{v_j,v_k\}\in EG} \Delta_{jk}^4 + \frac{R(\vec{\xi})}{2}
			\Bigg) 
		 d\vec{\xi},\\
		  	R(\vec{\xi}) = \mbox{tr} (\Lambda(\vec{\xi})  \hat{Q}^{-1} \Lambda(\vec{\xi})  \hat{Q}^{-1}),\\
		   \ \ \vec{\beta} = Q\vec{\alpha},
		\end{aligned} 
	\end{equation}
where $\Lambda(\vec{\xi})$ denotes the diagonal
matrix whose diagonal elements are equal to corresponding  components of the vector $Q\vec{\xi}$ and
  $\vec{\alpha}$ is the vector composed of the diagonal elements of $\hat{Q}^{-1}$.

Let define $\vec{\xi}(\vec{\theta}) = (\xi_1(\vec{\theta}), \xi_2(\vec{\theta}), \ldots \xi_n(\vec{\theta}))$ by
\begin{equation}\label{eq5.6}
	\theta_k = \sqrt{(d_k+1)/2} \,\xi_k.
\end{equation}
Then we can rewrite \eqref{eq5.5} in notations of Section 4:
\begin{equation}\label{eq5.7}
	\mbox{Int} = <e^{i\vec{p}^T\vec{\theta}}>_{F,\Omega} \prod\limits_{j=1}\limits^{n} \frac{1}{\sqrt{(d_j+1)/2}}, 
\end{equation}
where $\vec{p} = (p_1, p_2, \ldots p_n)$,
\begin{equation}\label{eq5.8}
	p_k = \frac{\beta_k}{\sqrt{(d_j+1)/2}},
\end{equation}
\begin{equation}
	F(\vec{\theta}) = - \vec{\theta}^T A \vec{\theta} + H(\vec{\theta}),
\end{equation}
\begin{equation}\label{eq5.10}
	\vec{\theta}^T A \vec{\theta}= \frac{1}{2}\vec{\xi}(\vec{\theta})^T \hat{Q}\vec{\xi}(\vec{\theta}) , 
	\ \ \ \ A_{jk} = \frac{1}{\sqrt{(d_j+1)(d_k+1)}}\hat{Q}_{jk},
\end{equation}
\begin{equation}\label{eq5.11}
	H(\theta) = 	- \frac{1}{12}\sum\limits_{\{v_j,v_k\}\in EG} \Delta_{jk}^4 + \frac{R(\vec{\xi}(\vec{\theta}))}{2},
\end{equation}
\begin{equation}\label{eq5.12}
	\Omega = \{\vec{\theta} \in \mathbb{R}^n : \vec{\xi}(\vec{\theta}) \in U_n(n^{-1/2+\varepsilon})\}.
\end{equation} 		

We aim to reduce, using Lemma \ref{Proposition_4.1}, expression $\mbox{Int}$ of \eqref{eq5.7} to 
\begin{equation}\label{eq5.13}
	<1>=   \int\limits_{\mathbb{R}^n}  e^{- \vec{\theta}^T A \vec{\theta}} d\vec{\theta} = \frac{\pi^{n/2}}{\sqrt{\det A}} = 
	\frac{(2\pi)^{n/2}}{\sqrt{\det \hat{Q}}} \prod\limits_{j=1}\limits^{n}{\sqrt{(d_j+1)/2}}. 	
\end{equation}
Our argument is as follows: first we have to verify that 
all assumptions of Lemma \ref{Proposition_4.1} hold, then we will gradually get rid of the oscillating term $e^{i\vec{p}^T\vec{\theta}}$, 
quadratic term $\frac{R(\vec{\xi}(\vec{\theta}))}{2}$ and the residual term $- \frac{1}{12}\sum\limits_{\{v_j,v_k\}\in EG} \Delta_{jk}^4$.

Further, we always use notation $f = O(g)$ meaning that $|f|\leq c|g|$ for some 
constant $c>0$ depending only on $\gamma$ and $\varepsilon$.

%%%%%%%%%%%%%%%%%%%%%%%%%%%%%%%%%%%%%%%%%%%%%%%%%%%%%%%%%%%%%%%%%%%%%%%%%%%%%%%%%%%%%%%%%%%%%%%%%%%%%%%%%%%%%%%%%%%%%%%%%
%%%%%%%%%%%%%%%%%%%%%%%%%%%%%%%%%%%%%%%%%%%%%%%%%%%%%%%%%%%%%%%%%%%%%%%%%%%%%%%%%%%%%%%%%%%%%%%%%%%%%%%%%%%%%%%%%%%%%%%%%

\subsection{Assumptions of Lemma \ref{Proposition_4.1}}

Combining  \eqref{4.2}, \eqref{4.3}, \eqref{eq5.6}, \eqref{eq5.10} and \eqref{eq5.12}, we get that
\begin{equation}\label{eq5.14}
\text{$A$ satisfy (\ref{A_ass}) } \ \text{ and } \
	U_n(r_1n^{\varepsilon}) \subset \Omega \subset U_n(r_2n^{\varepsilon}) 
\end{equation} 
for some $a,b,r_1,r_2>0$ depending only on $\gamma$.

Let $\vec{e}^{(k)} = (e^{(k)}_1,\ldots,e^{(k)}_n)  \in \mathbb{R}^n$ be defined by $e^{(k)}_j = \delta_{jk}$, where
$\delta_{jk}$ is the Kronecker delta. Due to the linearity of $\Lambda(\vec{\xi})$ and $\mbox{tr}(\cdot)$, we find that 
\begin{equation}\label{eq5.15}
 R(\vec{\xi}(\vec{\theta})) = \vec{\xi}(\vec{\theta})^T R \vec{\xi}(\vec{\theta}) = \vec{\theta}^T S \vec{\theta}, 
\end{equation}
where 
\begin{equation}\label{eq5.16}
	\begin{aligned}
		&R_{jk} = \mbox{tr} (\Lambda(\vec{e}^{(j)} )  \hat{Q}^{-1} \Lambda(\vec{e}^{(k)})  \hat{Q}^{-1}),\\
		&S_{jk} = \frac{R_{jk}}{{\sqrt{(d_j+1)/2}}{\sqrt{(d_k+1)/2}}}.
	\end{aligned}	
\end{equation}
We use the following inequalities for $n\times n$ matrices $X,Y$:
\begin{equation}\label{eq5.17}
	\begin{aligned}
		|\mbox{tr}(XY)| \leq \|X\|_{HS}\|Y\|_{HS},\\
			\|XY\|_{HS} \leq \|X\|_{HS}\|Y^T\|_{2},
	\end{aligned}
\end{equation}
where  $\|\cdot\|_{HS}$ denotes the Hilbert-Schmidt norm,
\begin{equation}
	\|X\|_{HS} = \sqrt{\sum\limits_{j=1}\limits^{n}\sum\limits_{k=1}\limits^{n} |X_{jk}|^2}.
\end{equation}
Combining  \eqref{eq5.3}, \eqref{4.2}, \eqref{eq5.15}-\eqref{eq5.17}, we find that
\begin{equation}\label{eqR_jk}
	\begin{aligned}
	\vec{\xi}_1^T R \vec{\xi}_2	 &\leq \|\Lambda(\vec{\xi_1})  \hat{Q}^{-1}\|_{HS} \|\Lambda(\vec{\xi_2})  \hat{Q}^{-1}\|_{HS} \leq
	\\ &\leq \|\Lambda(\vec{\xi}_1)\|_{HS}\|\Lambda(\vec{\xi}_2)\|_{HS}\|\hat{Q}^{-1}\|_2^2 =
	 \|Q\vec{\xi}_1\|_2\|Q\vec{\xi}_2\|_2\|\hat{Q}^{-1}\|_2^2 \leq\\ &\leq \|Q\|_2^2\|\hat{Q}^{-1}\|_2^2 \|\vec{\xi}_1\|_2\|\vec{\xi}_2\|_2 
	= O(1) \|\vec{\xi}_1\|_2\|\vec{\xi}_2\|_2.
	\end{aligned}
\end{equation}
Using \eqref{eq5.3}, \eqref{eq5.10}, \eqref{eq5.11} and \eqref{eqR_jk}, we obtain that
\begin{equation}\label{eq5.19}
	\begin{aligned}
	H(\vec{\theta}) \leq R(\vec{\xi}(\vec{\theta})) =\vec{\xi}(\vec{\theta})^T R \vec{\xi}(\vec{\theta}) = 	 
	O(1) \|\vec{\xi}(\vec{\theta})\|_2^2 = \\ =
	O(1) \frac{\vec{\xi}(\vec{\theta})^T\hat{Q}\vec{\xi}(\vec{\theta})}{n} = O(n^{-1}) \vec{\theta}^T A \vec{\theta}.
	\end{aligned}
\end{equation}

Let $(\hat{Q}^{-1}\Lambda(\vec{\xi}) \hat{Q}^{-1})_{kk}$ denote the $(k,k)$-th element of the matrix $\hat{Q}^{-1}\Lambda(\vec{\xi}) \hat{Q}^{-1}$. 
For any $1\leq k\leq n$, we have that
 \begin{equation}\label{eq5.20}
 	\begin{aligned}
 		\frac{\partial R(\vec{\xi})}{\partial \xi_k} &= 
 		2 \mbox{tr}\left(\frac{\partial \Lambda(\vec{\xi})}{ \partial \xi_k} \hat{Q}^{-1}\Lambda(\vec{\xi}) \hat{Q}^{-1}\right) =
 		2 \mbox{tr}\left( \Lambda(\vec{e}^{(k)}) \hat{Q}^{-1}\Lambda(\vec{\xi}) \hat{Q}^{-1}\right)=\\
 		&= 2 d_k (\hat{Q}^{-1}\Lambda(\vec{\xi}) \hat{Q}^{-1})_{kk} + 2\mbox{tr} 
 		\left(\tilde{\Lambda} \hat{Q}^{-1}\Lambda(\vec{\xi}) \hat{Q}^{-1}\right),
 		\end{aligned}
 \end{equation}
 where $\tilde{\Lambda}$ is the diagonal matrix with the diagonal elements $\tilde{\Lambda}_{jj} = \Lambda(\vec{e}^{(k)})_{jj}$ 
 in the case of $j\neq k$ and   $\tilde{\Lambda}_{kk} =0$. In particular, we have that
 \begin{equation}\label{eq5.21}
 	\|\tilde{\Lambda}\|_2 \leq 1.
 \end{equation}
Since $\Lambda(\vec{\xi})$ is diagonal matrix,  we get that
\begin{equation}
	|d_k (\hat{Q}^{-1}\Lambda(\vec{\xi}) \hat{Q}^{-1})_{kk}| = |d_k \Lambda(\vec{\xi})_{kk}| \|(\hat{Q}^{-1})_k\|_2^2 \leq 
	d_k \|\Lambda(\vec{\xi})\|_2 \|(\hat{Q}^{-1})_k\|_2^2,
\end{equation}
 where $(\hat{Q}^{-1})_k$ is the $k$-th column of the matrix $\hat{Q}^{-1}$. Note that
\begin{equation}
	\|(\hat{Q}^{-1})_k\|_2 \leq \|\hat{Q}^{-1}\|_2 \|\hat{Q}(\hat{Q}^{-1})_k\|_2  = \|\hat{Q}^{-1}\|_2.
\end{equation}
We also note that 
\begin{equation}\label{eq5.24}
	\|\Lambda(\vec{\xi})\|_2 = \|Q\vec{\xi}\|_\infty \leq 2n \|\vec{\xi}\|_{\infty}.
\end{equation}
Combining \eqref{4.2}, \eqref{4.3}, \eqref{eq5.21}-\eqref{eq5.24}, we get that
\begin{equation}\label{eq5.25}
	\begin{aligned}
		|\mbox{tr} 
 		\left(\tilde{\Lambda} \hat{Q}^{-1}\Lambda(\vec{\xi}) \hat{Q}^{-1}\right)| \leq
 		n\|\tilde{\Lambda} \hat{Q}^{-1}\Lambda(\vec{\xi}) \hat{Q}^{-1}\|_2 \leq \\ 
 		\leq n \|\tilde{\Lambda}\|_2 \|\Lambda(\vec{\xi})\|_2 \|\hat{Q}^{-1}\|_2^2 = O(1)\|\vec{\xi}\|_{\infty},\\
 				d_k (\hat{Q}^{-1}\Lambda(\vec{\xi}) \hat{Q}^{-1})_{kk} = O(1)\|\vec{\xi}\|_{\infty}.
	\end{aligned}
\end{equation}
Using \eqref{4.3}, \eqref{eq5.6}, \eqref{eq5.15}, \eqref{eq5.20} and \eqref{eq5.25} for all $1\leq k\leq n$, we obtain that
\begin{equation}\label{eq5.26}
	\begin{aligned}
		2\|R\vec{\xi}\|_\infty = \left\|\frac{\partial R(\vec{\xi})} {\partial \vec{\xi}}\right\|_{\infty} = 
		 O(1) \|\vec{\xi}\|_\infty,\\
			2\|S\vec{\theta}\|_\infty = 
			\left\|\frac{\partial R(\vec{\xi}(\vec{\theta}))} {\partial \vec{\theta}}\right\|_{\infty} = O(n^{-1})  \|\vec{\theta}\|_\infty.
	\end{aligned}
\end{equation}

For any $1\leq k\leq n$, we have that
\begin{equation}\label{eq5.27}
		\left|\frac{\partial}{\partial \xi_k}  \sum\limits_{\{v_j,v_l\}\in EG} \Delta_{jl}^4\right| \leq
		4\sum\limits_{j=1}\limits^{n}  |\Delta_{jk}|^3 = O(n) \|\vec{\xi}\|_\infty^3
\end{equation}
Combining \eqref{eq5.6} and \eqref{eq5.27}, we find that
\begin{equation}\label{eq5.28}
	\left\|\frac{\partial}{\partial \vec{\theta}}  \sum\limits_{\{v_j,v_l\}\in EG} \Delta_{jl}^4\right\|_\infty = O(n^{-1})\|\vec{\theta}\|_\infty^3.
\end{equation}

	Using \eqref{3.10}, the fact that
\begin{equation}\label{5.29}
	(\hat{Q}^{-1})_{jk} = \frac{(A^{-1})_{jk}}{{\sqrt{(d_j+1)/2}}{\sqrt{(d_k+1)/2}}},
\end{equation}
and \eqref{4.3}, we get that 
\begin{equation}\label{eq5.30}
\|\vec{p}\|_\infty = \sup\limits_{1\leq k\leq n} \frac{\beta_k}{\sqrt{(d_j+1)/2}} = O(n^{-1/2}) \|Q \alpha\|_\infty = O(n^{-1/2}),
\end{equation}
where vector $\vec{p}$ is the same as in \eqref{eq5.8}.

Putting together \eqref{eq5.10}, \eqref{eq5.14}, \eqref{eq5.19}, \eqref{eq5.26}, \eqref{eq5.28} and \eqref{eq5.30}, we get that
all assumptions of Lemma \ref{Proposition_4.1} with data \eqref{eq5.8}- \eqref{eq5.12}
hold for some constants $a, b, r_1, r_2, c_1, c_2 >0$ depending only on $\gamma$.	
		
%%%%%%%%%%%%%%%%%%%%%%%%%%%%%%%%%%%%%%%%%%%%%%%%%%%%%%%%%%%%%%%%%%%%%%%%%%%%%%%%%%%%%%%%%%%%%%%%%%%%%%%%%%%%%%%%%%%%%%%%%
%%%%%%%%%%%%%%%%%%%%%%%%%%%%%%%%%%%%%%%%%%%%%%%%%%%%%%%%%%%%%%%%%%%%%%%%%%%%%%%%%%%%%%%%%%%%%%%%%%%%%%%%%%%%%%%%%%%%%%%%%

\subsection{Oscillating term}		
Let define $\vec{p}^{(k)} = (p^{(k)}_1, p^{(k)}_2, \ldots, p^{(k)}_n) \in \mathbb{R}^n$ as follows:
\begin{equation}
	p^{(k)}_j = \left\{ 
	\begin{array}{lr}
	 		 	0,	& \text{ for } j \leq k, \\
	 		 	p_j, & \text{  otherwise}. 
	\end{array}
	\right.
\end{equation}
Let
\begin{equation}
	\begin{aligned}
	F^{(k)}(\vec{\theta}) &= 	 F(\vec{\theta}) 		-\frac{1}{2}\sum\limits_{j=1}\limits^{k} p_j^2\theta_j^2 .
	\end{aligned}
\end{equation}
Note that
\begin{equation}\label{F0=F}
	F^{(0)} (\vec{\theta}) = F(\vec{\theta}).
\end{equation}

Using \eqref{eq5.30} and  the fact that
	\begin{equation}
		\frac{\partial}{\partial \theta_k}\sum\limits_{j=1}\limits^{n}  p_j^2\theta_j^2  = 2 p_k^2 \theta_k = O(n^{-1}) \|\theta\|_\infty,
	\end{equation}
 we find that all assumptions of Lemma \ref{Proposition_4.1} hold  for 
	the case when we take function $F^{(k)}(\vec{\theta})$ instead of the function $F(\vec{\theta})$, vector $\vec{p}^{(k)}$ instead of vector $\vec{p}$
 and	for some constants $a, b, r_1, r_2, c_1, c_2 >0$ depending only on $\gamma$.

Note that 
\begin{equation}\label{eq5.33}
  \theta \in \Omega \ \Longrightarrow \ \|\theta\|_\infty = O(n^{\varepsilon}).
\end{equation}  
  Taking into account \eqref{eq5.30} and using the Taylor series expansion, we get that
\begin{equation}\label{eq5.34}
	e^{ip_k\theta_k} = e^{-\frac{1}{2} p_k^2 \theta_k^2} + i p_k\theta_k + O(n^{-3/2 + 3\varepsilon}), \ \ \vec{\theta}\in \Omega.
\end{equation}
Using \eqref{eq5.30} and \eqref{theta^0}, \eqref{theta+o} for $F^{(k-1)}$ and $\vec{p}^{(k-1)}$, we find that
\begin{equation}\label{eq5.35}
	\begin{aligned}
	i p_k <\theta_k e^{i \vec{\theta}^T\vec{p}^{(k)} }&>_{F^{{(k-1)}},\Omega} =\\ 
	= - \frac{p_k}{2}  \sum\limits_{j=k+1}\limits^n p_j &(A^{-1})_{jk} 
	 		< e^{i \vec{\theta}^T\vec{p}^{(k)}}>_{F^{(k)}, \Omega} 
	 		+ O(n^{-3/2+5\varepsilon})<1>,
	\end{aligned}
\end{equation}
 Combining \eqref{theta^0}, \eqref{eq5.34} and \eqref{eq5.35}, we obtain that
\begin{equation}\label{eq5.36}
	\begin{aligned}
	<e^{i \vec{\theta}^T\vec{p}^{(k-1)}}>_{F^{(k-1)},\Omega} &= \\=
	\Bigg(1 -  \frac{p_k}{2}  \sum\limits_{j=k+1}\limits^n p_j (A^{-1})_{jk} \Bigg) <&e^{i \vec{\theta}^T\vec{p}^{(k)}}>_{F^{(k)},\Omega} +\\
	&+ O(n^{-3/2+5\varepsilon})<1>.
	\end{aligned}
\end{equation}
Using \eqref{3.10}, \eqref{eq5.8}, \eqref{eq5.10} and \eqref{eq5.30}, we note also that
\begin{equation}\label{eq5.37}
	\begin{aligned}
	p_j p_k (A^{-1})_{jk} = \beta_j \beta_k (\hat{Q}^{-1})_{jk} = O(n^{-2}), \ \ j\neq k, \\
		p_k  \sum\limits_{j\neq k, j\leq n} p_j (A^{-1})_{jk}  = O(n^{-1}).
	\end{aligned}
\end{equation}
%We gradually get rid of the oscillations:
%	\begin{equation}
%		\vec{\theta}^T\vec{p} \rightarrow (\vec{\theta}^T\vec{p}- {p_1 \theta_1})
%		\rightarrow (\vec{\theta}^T\vec{p}- {p_1 \theta_1} - {p_2 \theta_2}) \rightarrow \ldots \rightarrow 0.
%	\end{equation} 
%	At this time, function $F$ is modified as follows:
%		\begin{equation}\label{eq5.39}
%		F \rightarrow ( F - \frac{1}{2}{p_1^2 \theta_1^2})
%		 \rightarrow \ldots \rightarrow 
%		(F - \frac{1}{2}p_1^2 \theta_1^2 - \ldots - \frac{1}{2} p_n^2 \theta_n^2).
%	\end{equation}
		Using  \eqref{F0=F}, \eqref{eq5.37} and \eqref{eq5.36} for $k=1,2\ldots n$,  we get that
	\begin{equation}\label{eq5.41}
		\begin{aligned}
				<e^{i \vec{\theta}^T\vec{p}}>_{F,\Omega} = 
	C_1	<1>_{F^{(n)},\Omega} 
	+ O(n^{-1/2+5\varepsilon})<1>,
		\end{aligned}
	\end{equation}
	where 
	\begin{equation}\label{C_1+++}
		\begin{aligned}
		C_1 &= \exp\left(- \sum_{k=1}^{n-1}\sum\limits_{j=k+1}\limits^{n}  {\beta_k}  (\hat{Q}^{-1})_{jk}  \beta_j\right)= \\&=
		\prod_{k=1}^{n-1}\left(1 - \frac{p_k}{2} \sum\limits_{j=k+1}\limits^{n} (\hat{A}^{-1})_{jk} p_j\right) 
		+ O(n^{-1})=  O(1).
		\end{aligned}
	\end{equation}
	
	  Taking into account \eqref{eq5.30} and using the Taylor series expansion, we get that
\begin{equation}\label{eq5.43}
	e^{-\frac{1}{2} p_k^2 \theta_k^2} =  1 - \frac{1}{2} p_k^2 \theta_k^2 + O(n^{-2 + 4\varepsilon}), \ \ \vec{\theta}\in \Omega.
\end{equation}
	Combining  \eqref{eq5.43} and \eqref{theta^0}, \eqref{theta^2} for $F^{{(k)}}$,  we get that
\begin{equation}\label{eq5.44}
	<1>_{F^{(k)},\Omega} = \left(1-\frac{1}{4}p_k^2\right)<1>_{F^{{(k-1)}},\Omega}+O(n^{-2 + 4\varepsilon})<1>.
\end{equation}
Using  \eqref{eq5.30}, \eqref{F0=F} and  \eqref{eq5.44} for $k=1,2\ldots,n$, we find that
\begin{equation}\label{eq5.45}
	<1>_{F^{(n)},\Omega}  = 
	C_2 <1>_{F,\Omega} + O(n^{-1+4\varepsilon})<1>, 
	\end{equation}
where 
\begin{equation}\label{C_2+++}
	C_2 = \exp\left(-\sum_{k=1}^{n} \frac{\beta_k^2}{2(d_k+1)}\right) =\prod_{k=1}^{n} \left(1 - \frac{1}{4}p_k^2\right) + O(n^{-1})= O(1).
\end{equation}	%%%%%%%%%%%%%%%%%%%%%%%%%%%%%%%%%%%%%%%%%%%%%%%%%%%%%%%%%%%%%%%%%%%%%%%%%%%%%%%%%%%%%%%%%%%%%%%%%%%%%%%%%%%%%%%%%%%%%%%%%
%%%%%%%%%%%%%%%%%%%%%%%%%%%%%%%%%%%%%%%%%%%%%%%%%%%%%%%%%%%%%%%%%%%%%%%%%%%%%%%%%%%%%%%%%%%%%%%%%%%%%%%%%%%%%%%%%%%%%%%%%
\subsection{Quadratic term}
Let define $\vec{\theta}^{k} = (\theta^{k}_1, \theta^{k}_1, \ldots, \theta^{k}_n) \in \mathbb{R}^n$ as follows:
\begin{equation}
	\theta^{k}_j = \left\{ 
	\begin{array}{lr}
	 		 	0,	& \text{ for } j \leq k, \\
	 		 	\theta_j, & \text{  otherwise}. 
	\end{array}
	\right.
\end{equation}
Let
\begin{equation}\label{eq5.48}
	\begin{aligned}
	F^{k}(\vec{\theta}) &= 	- \vec{\theta}^T A \vec{\theta} + H^{k}(\vec{\theta}),\\
		H^{k}(\vec{\theta}) &= - \frac{1}{12}\sum\limits_{\{v_j,v_k\}\in EG} \Delta_{jk}^4 + \frac{1}{2}R(\vec{\xi}(\vec{\theta}^{k})).
	\end{aligned}
\end{equation}
In absolutely similar way as given in Subsection 5.1, we find that all assumptions of Lemma \ref{Proposition_4.1} hold  for 
	the case when we take $F^{k}(\vec{\theta})$ instead of $F(\vec{\theta})$ for some constants $a, b, r_1, r_2, c_1, c_2 >0$ depending only on $\gamma$. 

Note that
	\begin{equation}\label{F00=F}
		F^{0}(\vec{\theta}) = F(\vec{\theta}),
	\end{equation}
\begin{equation}
	R(\vec{\xi}(\vec{\theta}^{k})) = \vec{\xi}(\vec{\theta}^{k})^T R \vec{\xi}(\vec{\theta}^{k}) 
	= {\vec{\theta}^{kT}}  S \vec{\theta}^{k},
\end{equation}
where matrices $R$, $S$ are the same that in \eqref{eq5.15}, \eqref{eq5.16}.

Combining \eqref{eq5.48}, the Mean Value Theorem and  \eqref{eq5.26}, we get that
\begin{equation}\label{eq5.52}
	\begin{aligned}
	|F^{k-1}(\vec{\theta}) - F^{k}(\vec{\theta})| = \frac{1}{2} | R(\vec{\xi}(\vec{\theta}^{k-1})) - R(\vec{\xi}(\vec{\theta}^{k}))| =\\
	=\frac{1}{2}\left|\frac{\partial R(\vec{\xi}(\vec{\theta}^*))} {\partial \theta_k} \theta_k\right| =
	 O(n^{-1}) \|\vec{\theta}\|_\infty^2,
	\end{aligned}
\end{equation}
where $\vec{\theta}^*$ lies on the segment between $\vec{\theta}^{k-1}$ and $\vec{\theta}^{k}$. 
Using \eqref{eq5.33}, \eqref{eq5.52}  and the Taylor series expansion, we get that
\begin{equation}\label{eq5.53}
	\begin{aligned}
	e^{F^{k-1}(\vec{\theta}) - F^{k}(\vec{\theta})} = 1 + F^{k-1}(\vec{\theta}) - F^{k}(\vec{\theta}) + O(n^{-2 + 4\varepsilon}) =\\
	= 1 + \frac{1}{2}S_{kk} \theta_k^2 + \sum\limits_{j=k+1}\limits^n S_{jk} \theta_j\theta_k + O(n^{-2 + 4\varepsilon}), 
	\ \ \vec{\theta}\in \Omega.
	\end{aligned}
\end{equation}
Using \eqref{4.3}, \eqref{eq5.16}, \eqref{eqR_jk}, we find that
\begin{equation}\label{eq5.51}
	S_{kk} = \frac{\vec{e}^{(k)T} R \vec{e}^{(k)}}{{(d_k+1)/2}} = O(n^{-1})\|\vec{e}^{(k)}\|_2^2 = O(n^{-1}).
\end{equation}
Since \eqref{eq5.26} imply that $\|S\|_\infty = O(n^{-1})$, using \eqref{theta^11}, we get that
\begin{equation}\label{eq5.54}
	\begin{aligned}
	\left|\sum\limits_{j=k+1}\limits^n S_{jk} <\theta_j\theta_k>_{F^{k-1}, \Omega} \right| \leq 
	 |\sup\limits_{k<j\leq n} <\theta_j\theta_k>_{F^{k-1}, \Omega}| \sum\limits_{j=k+1}\limits^n |S_{jk}| = \\
	 = O(n^{-1+5\varepsilon})\|S\|_\infty = O(n^{-2+5\varepsilon}).
	\end{aligned}
\end{equation}
Combining \eqref{theta^0}, \eqref{theta^2}, \eqref{eq5.53} and \eqref{eq5.54}, we obtain that
\begin{equation}\label{eq5.54+}
	<1>_{F^{k-1}, \Omega} = \left(1 + \frac{1}{4}S_{kk}\right)<1>_{F^{k},\Omega} + O(n^{-2+5\varepsilon})<1>.
\end{equation}
Using  \eqref{F00=F}, \eqref{eq5.51} and  \eqref{eq5.54+} for $k=1,2\ldots,n$, we find that
\begin{equation}\label{eq5.56}
	<1>_{F,\Omega}  = 
	C_3 <1>_{F^{n},\Omega} + O(n^{-1+5\varepsilon})<1>, 
	\end{equation}
where 
\begin{equation}\label{C_3+++}
		C_3  = \exp\left( \sum_{k=1}^{n}  \frac{R_{kk}}{2(d_k+1)}\right) = \prod_{k=1}^{n} \left(1 + \frac{1}{4}S_{kk}\right) + O(n^{-1})= O(1).
\end{equation}	%%%%%%%%%%%%%%%%%%%%%%%%%%%%%%%%%%%%%%%%%%%%%%%%%%%%%%%%%%%%%%%%%%%%%%%%%%%%%%%%%%%%%%%%%%%%%%%%%%%%%%%%%%%%%%%%%%%%%%%%%
%%%%%%%%%%%%%%%%%%%%%%%%%%%%%%%%%%%%%%%%%%%%%%%%%%%%%%%%%%%%%%%%%%%%%%%%%%%%%%%%%%%%%%%%%%%%%%%%%%%%%%%%%%%%%%%%%%%%%%%%%

\subsection{Residual term}
For a subset  $\Theta$ of $EG$ we define 
\begin{equation}
	\begin{aligned}
	F_\Theta(\vec{\theta}) &= 	- \vec{\theta}^T A \vec{\theta} + H_\Theta(\vec{\theta}),\\
		H_\Theta(\vec{\theta}) &= - \frac{1}{12}\sum\limits_{\{v_j,v_k\}\in \Theta} \Delta_{jk}^4.
	\end{aligned}
\end{equation}
In absolutely similar way as given in Subsection 5.1, we find that all assumptions of Lemma \ref{Proposition_4.1} hold  for 
	the case when we take $F_\Theta(\vec{\theta})$ instead of $F(\vec{\theta})$ for some constants $a, b, r_1, r_2, c_1, c_2 >0$ depending only on $\gamma$. 

Note that
\begin{equation}\label{Fn=0}
	F_{EG}(\vec{\theta}) = F^{n}(\vec{\theta}),
\end{equation}
\begin{equation}\label{4.7}
		\begin{aligned}
			\Delta_{jk}^4 = \frac{4\theta_j^4}{(d_j+1)^2}
			- 4 \frac{4\theta_j^3\theta_k}{(d_j+1)^{3/2}(d_k+1)^{1/2}}
			+ 6 \frac{4\theta_j^2\theta_k^2}{(d_j+1)(d_k+1)}-\\
			- 4 \frac{4\theta_k^3\theta_j}{(d_k+1)^{3/2}(d_j+1)^{1/2}}
			+ \frac{4\theta_k^4}{(d_k+1)^2},
			\end{aligned}
\end{equation}
Combining (\ref{theta^4}), (\ref{theta^13}), (\ref{theta^22}), \eqref{4.3} and (\ref{4.7}), we get that
	\begin{equation}\label{eq5.61}
			\begin{aligned}
			<\Delta_{jk}^4>_{F_{\Theta},\Omega} =\Bigg(\frac{3}{4}&\frac{4}{(d_j+1)^2} + \frac{6}{4}\frac{4}{(d_j+1)(d_k+1)} +\\&+ 	
			\frac{3}{4}\frac{4}{(d_k+1)^2}\Bigg) <1>_{F_{\Theta},\Omega} +
			  O(n^{-3+7\varepsilon})<1>.
			\end{aligned}
	\end{equation}
Using (\ref{theta^0}) and \eqref{eq5.61}, we get that
\begin{equation}{\label{4.14}}
			\begin{aligned}
			<e^{-\frac{1}{12}\Delta_{jk}^4 }>_{F_{\Theta},\Omega} = 
				P_{jk} <1>_{F_{\Theta},\Omega} + O(n^{-3+7\varepsilon})<1>,
			\end{aligned}
\end{equation}
where
		\begin{equation}
			P_{jk} = 1 - \frac{1}{4(d_j+1)^2} -  \frac{1}{2(d_j+1)(d_k+1)} - \frac{1}{4(d_k+1)^2}.
		\end{equation}
Note that
\begin{equation}
	1-\frac{1}{n^2}\leq P_{jk}\leq 1.
\end{equation}
Using (\ref{theta^0}), (\ref{4.14}), we can  gradually remove all the edges from the residual term $H_{EG}$ and obtain that
\begin{equation}\label{4.17}
	<1>_{F_{EG}, \Omega} = \prod\limits_{\{v_j,v_k\}\in EG} P_{jk}  <1>_{\Omega} + O(n^{-1+7\varepsilon})<1>.
\end{equation}
Combining (\ref{Omega_infty}) and (\ref{4.17}), we get that
\begin{equation}\label{eq5.66}
	<1>_{F_{EG}, \Omega} =  C_4 <1> + O(n^{-1+7\varepsilon})<1>,
\end{equation}
where 
\begin{equation}\label{C_4+++}
C_4 = \exp\left(-\frac{1}{4}\sum\limits_{\{v_j,v_k\}\in EG} \left(\frac{1}{d_j+1} + \frac{1}{d_k+1}\right)^2 \right).
\end{equation}

Combining \eqref{eq5.41}, \eqref{eq5.45}, \eqref{eq5.56}, \eqref{Fn=0} and \eqref{eq5.66}, we find that
\begin{equation}\label{eq5.67}
	<e^{i\vec{p}^T\vec{\theta}}>_{F,\Omega} = C_1 C_2 C_3 C_4 <1> + O(n^{-1/2+7\varepsilon})<1>.
\end{equation}

Taking into account  \eqref{3.10}, \eqref{eq5.10}, \eqref{C_1+++}, \eqref{C_2+++}, we get that:
\begin{equation}
	\begin{aligned}
	C_1 C_2 = \exp\left(-\frac{1}{2}\vec{\beta}^T \hat{Q}^{-1} \vec{\beta} \right),\\
		\vec{\beta}^T \hat{Q}^{-1} \vec{\beta} = \vec{\alpha}^T Q^T  \hat{Q}^{-1} (Q+J)\vec{\alpha} 
		- \vec{\alpha}^T Q  \hat{Q}^{-1} J \vec{\alpha} = \\
		= \sum\limits_{\{v_j,v_k\}\in EG} \left(\frac{1}{d_j+1} - \frac{1}{d_k+1}\right)^2 + O(n^{-1}).
	\end{aligned}
\end{equation}
Using again \eqref{3.10}, \eqref{eq5.10}, we find that
\begin{equation}
	\begin{aligned}
		R_{kk} = \mbox{tr} (\Lambda(\vec{e}^{(k)} )  \hat{Q}^{-1} \Lambda(\vec{e}^{(k)})  \hat{Q}^{-1}) = 
		\\
		= \sum\limits_{j,m=1}\limits^n \Lambda_{jj}(\vec{e}^{(k)} )  (\hat{Q}^{-1})_{jm} \Lambda_{mm}(\vec{e}^{(k)})  
		(\hat{Q}^{-1})_{mj} 
		= 1 + O(n^{-1}),
	\end{aligned}
\end{equation} 
\begin{equation}\label{eq5.67fin}
	\begin{aligned}
	C_3 = \exp\left( \sum_{k=1}^{n}  \frac{R_{kk}}{2(d_k+1)}\right) 
	= \exp\left( \sum_{k=1}^{n}  \frac{1}{2(d_k+1)}\right) + O(n^{-1}) =\\
		=\exp\left( 
		 \frac{1}{2} \sum\limits_{\{v_j,v_k\}\in EG} \left( \frac{1}{(d_j+1)^2} + \frac{1}{(d_k+1)^2} \right)\right) + O(n^{-1}).
	\end{aligned}
\end{equation}

Putting together \eqref{EulS}, \eqref{S+S_0}, \eqref{S_0+++}, (\ref{tGQ}), \eqref{eq5.7}, \eqref{eq5.13} and \eqref{eq5.66}-\eqref{eq5.67fin},  we obtain (\ref{main_eq}) and 
\eqref{main_eq2} for $n\geq n_0(\gamma,\varepsilon)>0$ (with the exponent $7\varepsilon$ instead of $\varepsilon$).
Estimate \eqref{main_eq2} for $n\leq n_0$ can be fulfilled by choice of sufficiently large constant $C$.
	%%%%%%%%%%%%%%%%%%%%%%%%%%%%%%%%%%%%%%%%%%%%%%%%%%%%%%%%%%%%%%%%%%%%%%%%%%%%%%%%%%%%%%%%%%%%%%%%%%%%%%%%%%%%%%%%%%%%%%%%%
%%%%%%%%%%%%%%%%%%%%%%%%%%%%%%%%%%%%%%%%%%%%%%%%%%%%%%%%%%%%%%%%%%%%%%%%%%%%%%%%%%%%%%%%%%%%%%%%%%%%%%%%%%%%%%%%%%%%%%%%%
%%%%%%%%%%%%%%%%%%%%%%%%%%%%%%%%%%%%%%%%%%%%%%%%%%%%%%%%%%%%%%%%%%%%%%%%%%%%%%%%%%%%%%%%%%%%%%%%%%%%%%%%%%%%%%%%%%%%%%%%%
%%%%%%%%%%%%%%%%%%%%%%%%%%%%%%%%%%%%%%%%%%%%%%%%%%%%%%%%%%%%%%%%%%%%%%%%%%%%%%%%%%%%%%%%%%%%%%%%%%%%%%%%%%%%%%%%%%%%%%%%%
%%%%%%%%%%%%%%%%%%%%%%%%%%%%%%%%%%%%%%%%%%%%%%%%%%%%%%%%%%%%%%%%%%%%%%%%%%%%%%%%%%%%%%%%%%%%%%%%%%%%%%%%%%%%%%%%%%%%%%%%%	

\section{Proof of Lemma \ref{Proposition_4.1}}
In this section we use notation $f = O(g)$ meaning that $|f|\leq c|g|$ for some 
$c>0$ depending only on $r_1, r_2, c_1,c_2, a, b$ and $\varepsilon$.

Let 
\begin{equation}\label{5.1}
\vec{\phi}(\vec{\theta}) = (\phi_1(\vec{\theta}), \phi_2(\vec{\theta}), \ldots, \phi_n(\vec{\theta})) = A \vec{\theta}.
\end{equation}
According to (\ref{A_ass}), $A = I + X$,  $X_{jj}=0$,  and so
for some $g_1(\vec{\theta}) = g_1(\theta_2, \ldots, \theta_n)$ 
\begin{equation}\label{5.2}
	\vec{\theta}^T A \vec{\theta} = \phi_1^2(\vec{\theta}) + g_1(\vec{\theta}).
\end{equation} 
%Similarly, we define the functions $g_2, g_3, \ldots g_n$. 
Using (\ref{A_ass}), \eqref{5.2} and  estimating insignificant parts of Gaussian integral of the following type:
\begin{equation}
	\int \left(\max\{|x|,k_1\}\right)^s e^{-(x-k_2)^2} dx,
\end{equation}
we find that for $r>0$, $s\geq 0$
\begin{equation}\label{eq.6.3}
	\begin{aligned}
	<\|\vec{\theta}\|_\infty^s> &= \int\limits_{\mathbb{R}^n} \|\vec{\theta}\|_\infty^s e^{-\vec{\theta}^T A \vec{\theta}} d \vec{\theta}
	=\\&=
	\int\limits_{-\infty}\limits^{+\infty}\cdots \int\limits_{-\infty}\limits^{+\infty}
	 e^{-\,g_1(\theta_2,\ldots,\theta_{n})}
	\left( \int \limits_{-\infty}\limits^{+\infty} \|\vec{\theta}\|_\infty^s e^{-\phi_1(\vec{\theta})^2} d \theta_1  \right)
	d \theta_2 \ldots d \theta_n\\
	&=
		\left(1 + O\left(\exp(-c_4n^{2\varepsilon})\right)\right)
			\int \limits_{|\phi_1(\vec{\theta})|\leq r n^{\varepsilon}} \|\vec{\theta}\|_\infty^s e^{-\,\vec{\theta}^T A\vec{\theta}} d \vec{\theta},  
	\end{aligned}
\end{equation}
where $c_4 = c_4(r,\varepsilon,s)>0.$
Combining similar to \eqref{eq.6.3} expressions for $\phi_1, \phi_2, \ldots \phi_n$, we get that
\begin{equation}\label{5.4}
	\int \limits_{||\vec{\phi}(\vec{\theta})||_\infty\leq r n^{\varepsilon}} \|\vec{\theta}\|_\infty^s e^{-\vec{\theta}^T A \vec{\theta}} d \vec{\theta}
	= \left(1 + O\left(\exp(-c_5 n^{2\varepsilon})\right)\right) <\|\vec{\theta}\|_\infty^s>,
\end{equation}
where $c_5 = c_5(r,\varepsilon,s)>0.$ Combining (\ref{3.9}), (\ref{5.1}) and (\ref{5.4}) with $s=0$, we obtain (\ref{Omega_infty}).

Using (\ref{R_1}), we find that
\begin{equation}\label{5.5}
	|<1>_{F,\Omega}| \leq \int\limits_{\Omega} |e^{F(\vec{\theta})}| d\vec{\theta} 
	\leq 
	\int\limits_{\mathbb{R}^n} e^{-\vec{\theta}^T A \vec{\theta} + \frac{c_1}{n}\vec{\theta}^T A \vec{\theta}} d\vec{\theta} = O(<1>).
\end{equation}

In order to prove (\ref{theta^2}) - (\ref{theta^22}) we use the following two lemmas. The proofs of them are given in Section 7.
\begin{Lemma}\label{Lemma_5.1}
Let $\Omega\subset \mathbb{R}^n$ be such that  $U_n(r_1n^{\varepsilon}) \subset \Omega \subset U_n(r_2n^{\varepsilon})$ for some $r_2>r_1>0$. 
Let $A$ satisfy (\ref{A_ass}) and  assumptions (\ref{R_1}) and (\ref{R_2})  hold for some constants $c_1,c_2>0$. 
%Let 
% $T(\vec{\theta})$ be such that for any $\vec{\theta}\in{\mathbb{R}^n}$: 
%$|T(\vec{\theta})| \leq c_s n^s$, for some fixed constants $c_s, s > 0$. 
Let $P = P(x) = O(|x|^s)$ for some fixed $s\geq 0$. 
Then for any $T(\vec{\theta})$ such that $|T(\vec{\theta})| \leq P(\|\vec{\theta}\|_\infty):$ 
\begin{equation}
	\begin{aligned}
	<T(\vec{\theta})>_{\mathbb{R}^n \setminus \Omega} = O\left(\exp(-c_6n^{2\varepsilon})\right) <1>, 
	\end{aligned}
\end{equation}
and for any $\tilde{T}(\vec{\theta})=\tilde{T}(\theta_1,\ldots,\theta_{k-1},\theta_{k+1}, \ldots, \theta_n)$ such that
$\tilde{T}(\vec{\theta}) \leq P(\|\vec{\theta}\|_\infty)$:
\begin{equation}\label{phi^2}
	\begin{aligned}
	<\phi_k^2(\vec{\theta})\tilde{T}(\vec{\theta})>_{F, \Omega} = 
	\frac{1}{2}<\tilde{T}(\vec{\theta})>_{F,\Omega} + O(n^{-1+4\varepsilon})<|\tilde{T}(\vec{\theta})|>_{F,\Omega}+\\ 
	 + O\left(\exp(-c_6n^{2\varepsilon})\right)<1>,
	\end{aligned}
\end{equation}
\begin{equation}\label{phi^4}
	\begin{aligned}
	<\phi_k^4(\vec{\theta})\tilde{T}(\vec{\theta})>_{F, \Omega} = 
	\frac{3}{4}<\tilde{T}(\vec{\theta})>_{F,\Omega} + O(n^{-1+4\varepsilon})<|\tilde{T}(\vec{\theta})|>_{F,\Omega}+\\
	 + O\left(\exp(-c_6n^{2\varepsilon})\right)<1>,
	 \end{aligned}
\end{equation}
\begin{equation}\label{phi}
	\begin{aligned}
	<\phi_k(\vec{\theta})\tilde{T}(\vec{\theta})>_{F, \Omega} = 
	O(n^{-1+4\varepsilon}) <|\tilde{T}(\vec{\theta})|>_{F,\Omega} +\\
	 + O\left(\exp(-c_6n^{2\varepsilon})\right)<1>,
	 \end{aligned}
\end{equation}
\begin{equation}\label{phi^3}
	\begin{aligned}
	<\phi_k^3(\vec{\theta})\tilde{T}(\vec{\theta})>_{F, \Omega} = 
	O(n^{-1+6\varepsilon}) <|\tilde{T}(\vec{\theta})|>_{F,\Omega} +\\
	 + O\left(\exp(-c_6n^{2\varepsilon})\right)<1>,
	 \end{aligned}
\end{equation}
	where function $F$ is defined by (\ref{def_F}), 
	vector $\vec{\phi}(\vec{\theta})$ is defined by (\ref{5.1}) and constant $c_6 = c_6(r_1, r_2, c_1,c_2, a, b, \varepsilon, P)>0$.
\end{Lemma}
%Note that for any $r>0$
%\begin{equation}
%	\begin{aligned}
%	<1>_{r} = \int\limits_{U_n(rn^{\varepsilon})} \exp\left( - \vec{\theta}^T A \vec{\theta}\right) d \vec{\theta}= \
%	\end{aligned}
%\end{equation}
\begin{Lemma}\label{Lemma_5.2}
Let assumptions of Lemma \ref{Lemma_5.1} hold 
and $s_1,s_2,\ldots ,s_n \in \mathbb{N}\cup \{0\}$,
\begin{equation}
	\begin{aligned}
	M(\vec{x}) = x_1^{s_1}\cdots  x_{n}^{s_{n}},\\
	 s = s_1 + \ldots + s_n>0.
	\end{aligned}
\end{equation}
Let 
\begin{equation} 
	s_k =0 \text{ and  } |\{j : s_j \neq 0\}| \leq 3. 
\end{equation}
Then
\begin{equation}\label{5.17}
<\phi_k(\vec{\theta})M(\vec{\phi}(\vec{\theta}))>_{F, \Omega} = O(sn^{-1+(s+4)\varepsilon})<1>.
\end{equation}
\end{Lemma}

Using (\ref{3.10}), we find that
\begin{equation}\label{eq6.14}
	 \theta_k = \phi_k + \vec{q}_k^T \vec{\phi}, \ \ \|\vec{q}\|_{\infty} = O(n^{-1}).
\end{equation}
Combining  (\ref{5.5}), \eqref{eq6.14},  Lemma \ref{Lemma_5.1} and Lemma \ref{Lemma_5.2}, we obtain that:
\begin{equation}\label{5.15'}
\begin{aligned}
<&\delta_k(\vec{\theta})^2>_{F, \Omega} = <(\vec{q}_k^T \vec{\phi})^2>_{F, \Omega}=\\ &= O(n^{-2})
\left(
\sum\limits_{j} <\phi_j(\vec{\theta})^2>_{F, \Omega} +
\sum\limits_{j_1\neq j_2} |<\phi_{j_1}(\vec{\theta})\phi_{j_2}(\vec{\theta})>_{F, \Omega}|
\right) = \\&=
\left(O(n^{-1}) + O(n^{-1 + 5\varepsilon})\right)<1> = O(n^{-1 + 5\varepsilon})<1>,
\end{aligned}
\end{equation} 
\begin{equation}\label{5.16'}
\begin{aligned}
<\delta_k(\vec{\theta})^4>_{F, \Omega} &= <(\vec{q}_k^T \vec{\phi})^4>_{F, \Omega}= O(n^{-4})
\sum\limits_{j} <\phi_j(\vec{\theta})^4>_{F, \Omega} +  \\ &+
O(n^{-4})\sum\limits_{j_1\neq j_2} <\phi_{j_1}(\vec{\theta})^2\phi_{j_2}(\vec{\theta})^2>_{F, \Omega}
 + \\ &+ O(n^{-4})\sum\limits_{j_1\neq j_2} |<\phi_{j_1}(\vec{\theta})\phi_{j_2}(\vec{\theta})^3>_{F, \Omega}| +\\+
 O(n^{-4}&)\sum\limits_{j_1\neq j_2\neq j_3} |<\phi_{j_1}(\vec{\theta})\phi_{j_2}(\vec{\theta})\phi_{j_3}(\vec{\theta})^2>_{F, \Omega}| 
 +\\+ O(n^{-4})&\sum\limits_{j_1\neq j_2\neq j_3\neq j_4} 
 |<\phi_{j_1}(\vec{\theta})\phi_{j_2}(\vec{\theta})\phi_{j_3}(\vec{\theta})\phi_{j_4}(\vec{\theta})>_{F, \Omega}|
 =
  \\
  =
\Big(O(n^{-3}) + &O(n^{-2 + 4\varepsilon})+ O(n^{-3 + 7\varepsilon}) + O(n^{-2  + 7\varepsilon}) +\\&+ O(n^{-1 + 7\varepsilon})\Big) <1>
= O(n^{-1 + 7\varepsilon})<1>,
\end{aligned}
\end{equation} 
where 
\begin{equation}
	\delta_k(\vec{\theta}) = \theta_k - \phi_k(\vec{\theta}).
\end{equation}
According to (\ref{A_ass}), we have that
\begin{equation}\label{5.18'}
		\delta_k(\vec{\theta}) = 	\delta_k(\theta_1,\ldots,\theta_{k-1},\theta_{k+1}, \ldots, \theta_n) 
\end{equation}
Using (\ref{5.5}), \eqref{eq6.14}, (\ref{5.15'}), (\ref{5.16'}), (\ref{5.18'}) and Lemma \ref{Lemma_5.1}, we obtain that:
\begin{equation}\label{eq6.19}
	\begin{aligned}
	<\phi_k(\vec{\theta})\delta_k(\vec{\theta})>_{F, \Omega} &= \\=
	O(n^{-1+4\varepsilon})&<|\delta_k(\vec{\theta})|>_{F,\Omega}
	+ O\left(\exp(-c_6n^{2\varepsilon})\right)<1>=\\ &= O(n^{-1+5\epsilon})<1>,
	\end{aligned}
\end{equation}
\begin{equation}\label{eq6.20}
	\begin{aligned}
	<\theta_k^2>_{F, \Omega} = <(\phi_k(\vec{\theta})+&\delta_k(\vec{\theta}))^2>_{F, \Omega} = \\&=
	 <\phi_k(\vec{\theta})^2>_{F, \Omega} + O(n^{-1+5\epsilon})<1> = \\
	 &=\frac{1}{2}<1>_{F, \Omega}+O(n^{-1+5\epsilon})<1>.
	\end{aligned}
\end{equation}
\begin{equation}
	\begin{aligned}
	<\phi_k(\vec{\theta})\delta_k(\vec{\theta})^3>_{F, \Omega} &= \\=
	O(n^{-1+4\varepsilon})&<|\delta_k(\vec{\theta})|^3>_{F,\Omega}
	+ O\left(\exp(-c_6n^{2\varepsilon})\right)<1>=\\ &= O(n^{-1+7\epsilon})<1>,
	\end{aligned}
\end{equation}
\begin{equation}
	\begin{aligned}
	<\phi_k(\vec{\theta})^2\delta_k(\vec{\theta})^2>_{F, \Omega} = \
	\left(\frac{1}{2}+ O(n^{-1+4\varepsilon}) \right)<\delta_k(\vec{\theta})^2>_{F,\Omega}
	 +\\+ O\left(\exp(-c_6n^{2\varepsilon})\right)<1>=\\ = O(n^{-1+5\epsilon})<1>,
	\end{aligned}
\end{equation}
\begin{equation}\label{eq6.22}
	\begin{aligned}
	<\phi_k(\vec{\theta})^3\delta_k(\vec{\theta})>_{F, \Omega} &= \\
	=O(n^{-1+6\varepsilon})&<|\delta_k(\vec{\theta})|>_{F,\Omega}
	 + O\left(\exp(-c_6n^{2\varepsilon})\right)<1>=\\ &= O(n^{-1+7\epsilon})<1>,
	\end{aligned}
\end{equation}
\begin{equation}
	\begin{aligned}
	<\theta_k^4>_{F, \Omega} = <(\phi_k(\vec{\theta})+&\delta_k(\vec{\theta}))^4>_{F, \Omega} = \\&=
	 <\phi_k(\vec{\theta})^4>_{F, \Omega} + O(n^{-1+7\epsilon})<1> = \\
	 &=\frac{3}{4}<1>_{F, \Omega}+O(n^{-1+7\epsilon})<1>.
	\end{aligned}
\end{equation}
In a similar way as in (\ref{5.15'}), (\ref{5.16'}), using \eqref{eq6.19} - \eqref{eq6.22}, we find that
\begin{equation}\label{eq6.24}
\begin{aligned}
<\delta_k(\vec{\theta})\theta_l>_{F, \Omega} &= O(n^{-1})
\sum\limits_{j} <\phi_j(\vec{\theta})\theta_l>_{F, \Omega} =\\=
\Big(O(n^{-1+2\varepsilon}) &+ O(n^{-1 + 5\varepsilon})\Big)<1> = O(n^{-1 + 5\varepsilon})<1>,
\end{aligned}
\end{equation}
\begin{equation}\label{5.23}
\begin{aligned}
<\delta_k(\vec{\theta})^2\theta_l^2>_{F, \Omega} &= O(n^{-2})
\sum\limits_{j} <\phi_j(\vec{\theta})^2\theta_l^2>_{F, \Omega} +\\&+ O(n^{-2})
\sum\limits_{j_1\neq j_2} |<\phi_{j_1}(\vec{\theta})\phi_{j_2}(\vec{\theta})\theta_l^2>_{F, \Omega}| = \\=
\Big(O(n^{-1+2\varepsilon}) &+ O(n^{-1 + 7\varepsilon})\Big)<1> = O(n^{-1 + 7\varepsilon})<1>,
\end{aligned}
\end{equation}
\begin{equation}\label{5.24}
\begin{aligned}
<\delta_k(\vec{\theta})&\theta_l^3>_{F, \Omega} =\\&= O(n^{-1})
\sum\limits_{j\neq l} <\phi_j(\vec{\theta})\theta_l^3>_{F, \Omega} + O(n^{-1})<\phi_l(\vec{\theta})\theta_l^3>_{F, \Omega} = \\&=
\left(O(n^{-1+7\varepsilon}) + O(n^{-1 + 4\varepsilon})\right)<1> = O(n^{-1 + 7\varepsilon})<1>.
\end{aligned}
\end{equation}
Using (\ref{5.5}), (\ref{5.18'}), \eqref{eq6.24}-(\ref{5.24})  and Lemma \ref{Lemma_5.1}, we obtain that:
\begin{equation}
	\begin{aligned}
	<\phi_k(\vec{\theta})\delta_k(\vec{\theta})\theta_l^2>_{F, \Omega} = \
	O(n^{-1+4\varepsilon})<|\delta_k(\vec{\theta})\theta_l^2|>_{F,\Omega} + \\
	+ O\left(\exp(-c_6n^{2\varepsilon})\right)<1>=\\ = O(n^{-1+7\epsilon})<1>,
	\end{aligned}
\end{equation}
\begin{equation}
	\begin{aligned}
	<\theta_k\theta_l>_{F, \Omega} = <(\phi_k(\vec{\theta})+\delta_k(\vec{\theta}))\theta_l>_{F, \Omega} = 
	  O(n^{-1+5\epsilon})<1>,
	\end{aligned}
\end{equation}
\begin{equation}
	\begin{aligned}
	<\theta_k\theta_l^3>_{F, \Omega} = <(\phi_k(\vec{\theta})+\delta_k(\vec{\theta}))\theta_l^3>_{F, \Omega} = 
	  O(n^{-1+7\epsilon})<1>,
	\end{aligned}
\end{equation}
\begin{equation}
	\begin{aligned}
	<\theta_k^2\theta_l^2>_{F, \Omega} &= <(\phi_k(\vec{\theta})+\delta_k(\vec{\theta}))^2\theta_l^2>_{F, \Omega}\\ &= 
	 <\phi_k(\vec{\theta})^2\theta_l^2>_{F, \Omega} + O(n^{-1+7\epsilon})<1> =\\
	 &=\frac{1}{2}<\phi_l(\vec{\theta})+\delta_l(\vec{\theta}))^2>_{F, \Omega}+O(n^{-1+7\epsilon})<1> =\\&=\frac{1}{4}<1>_{F, \Omega}+O(n^{-1+7\epsilon})<1>.
	\end{aligned}
\end{equation}

Since $\|\vec{p}\|_\infty = O(n^{-1/2})$, using the fact that
\begin{equation}
	e^{\frac{1}{2}p_k^2\theta_k^2} = 1 + O(n^{-1+2\varepsilon}), \ \ \ \theta \in \Omega,
\end{equation}
and \eqref{5.5}, we get that
\begin{equation}\label{eq6.33}
	<\theta_k e^{i\vec{\theta}^T\vec{p} - i p_k \theta_k}>_{F,\Omega} = <\theta_k e^{i\vec{\theta}^T\vec{p} - i p_k \theta_k}>_{F',\Omega} + 
	O(n^{-1+3\varepsilon})<1>,
\end{equation}
where $F'=F - \frac{1}{2}p_k^2\theta_k^2$. It is clear that $F'$ satisfy all assumptions of Lemma \ref{Lemma_5.1}. 
For any $\vec{p}'= (p_1',p_2',\ldots, p_n')\in \mathbb{R}^n$ such that $\|\vec{p}'\|_\infty = O(n^{-1/2})$, we have that
\begin{equation}\label{eq6.34}
	e^{ip_l'\theta_l}  = 1 + ip_l'\theta_l+ O(n^{-1+2\varepsilon}), \ \ \, e^{-ip_l'\theta_l}= 1+O(n^{-1/2+\varepsilon}), \ \  \theta \in \Omega.
\end{equation}
Using  \eqref{5.5}, \eqref{phi^2} and \eqref{phi}  with $F'$ instead of $F$, \eqref{eq6.34}, we get that
\begin{equation}\label{eq6.35}
	\begin{aligned}
	<\phi_l(\vec{\theta}) e^{i\vec{\theta}^T\vec{p}'-ip_l'\theta_l}>_{F',\Omega} = O(n^{-1+4\varepsilon})<1>,\\
	<\phi_l(\vec{\theta}) \delta_l(\vec{\theta}) e^{i\vec{\theta}^T\vec{p}'-ip_l'\theta_l}>_{F',\Omega} = O(n^{-1+5\varepsilon})<1>,\\
	<\phi_l(\vec{\theta}) \theta_le^{i\vec{\theta}^T\vec{p}'-ip_l'\theta_l}>_{F',\Omega} = 
	<\phi_l^2(\vec{\theta}) e^{i\vec{\theta}^T\vec{p}'-ip_l'\theta_l}>_{F',\Omega} + \\
	<\phi_l(\vec{\theta}) \delta_l(\vec{\theta})e^{i\vec{\theta}^T\vec{p}'-ip_l'\theta_l}>_{F',\Omega}=\\
	= \frac{1}{2}<e^{i\vec{\theta}^T\vec{p}'-ip_l'\theta_l}>_{F',\Omega} + O(n^{-1+5\varepsilon})<1> =\\
	=\frac{1}{2}<e^{i\vec{\theta}^T\vec{p}'}>_{F',\Omega} + O(n^{-1/2+5\varepsilon})<1>.
	\end{aligned}
\end{equation}
Combining \eqref{eq6.34} and \eqref{eq6.35}, we find that
\begin{equation}\label{eq6.36}
	\begin{aligned}
 <\phi_l(\vec{\theta}) e^{i\vec{\theta}^T\vec{p}'}&>_{F',\Omega} = 
 <\phi_l(\vec{\theta}) e^{i\vec{\theta}^T\vec{p}'-ip_l'\theta_l}>_{F',\Omega} +\\
 &+ ip_l'<\phi_l(\vec{\theta}) \theta_le^{i\vec{\theta}^T\vec{p}'-ip_l'\theta_l}>_{F',\Omega} + O(n^{-1+3\varepsilon})<1> =\\
 &= \frac{i}{2}p_l'<e^{i\vec{\theta}^T\vec{p}'}>_{F',\Omega} + O(n^{-1+5\varepsilon})<1>.
 \end{aligned}
\end{equation}
Using \eqref{3.9}, \eqref{5.1} and \eqref{eq6.36}, we get that
\begin{equation}\label{eq6.37}
	\begin{aligned}
 			<\theta_k e^{i\vec{\theta}^T\vec{p}'}&>_{F',\Omega} = \\&=
 			\frac{i}{2} \sum\limits_{l=1}\limits^n p_l' (A^{-1})_{lk}<e^{i\vec{\theta}^T\vec{p}'}>_{F',\Omega} +
 			O(n^{-1+5\varepsilon})<1>
 \end{aligned}
\end{equation}
Combining \eqref{eq6.33} and \eqref{eq6.37} for $\vec{p}'= \vec{p} - p_k\vec{e}^{(k)}$, we obtain \eqref{theta+o}.
%%%%%%%%%%%%%%%%%%%%%%%%%%%%%%%%%%%%%%%%%%%%%%%%%%%%%%%%%%%%%%%%%%%%%%%%%%%%%%%%%%%%%%%%%%%%%%%%%%%%%%%%%%%%%%%%%%%%%%%%%
%%%%%%%%%%%%%%%%%%%%%%%%%%%%%%%%%%%%%%%%%%%%%%%%%%%%%%%%%%%%%%%%%%%%%%%%%%%%%%%%%%%%%%%%%%%%%%%%%%%%%%%%%%%%%%%%%%%%%%%%%
%%%%%%%%%%%%%%%%%%%%%%%%%%%%%%%%%%%%%%%%%%%%%%%%%%%%%%%%%%%%%%%%%%%%%%%%%%%%%%%%%%%%%%%%%%%%%%%%%%%%%%%%%%%%%%%%%%%%%%%%%
%%%%%%%%%%%%%%%%%%%%%%%%%%%%%%%%%%%%%%%%%%%%%%%%%%%%%%%%%%%%%%%%%%%%%%%%%%%%%%%%%%%%%%%%%%%%%%%%%%%%%%%%%%%%%%%%%%%%%%%%%
%%%%%%%%%%%%%%%%%%%%%%%%%%%%%%%%%%%%%%%%%%%%%%%%%%%%%%%%%%%%%%%%%%%%%%%%%%%%%%%%%%%%%%%%%%%%%%%%%%%%%%%%%%%%%%%%%%%%%%%%%
%%%%%%%%%%%%%%%%%%%%%%%%%%%%%%%%%%%%%%%%%%%%%%%%%%%%%%%%%%%%%%%%%%%%%%%%%%%%%%%%%%%%%%%%%%%%%%%%%%%%%%%%%%%%%%%%%%%%%%%%%
%%%%%%%%%%%%%%%%%%%%%%%%%%%%%%%%%%%%%%%%%%%%%%%%%%%%%%%%%%%%%%%%%%%%%%%%%%%%%%%%%%%%%%%%%%%%%%%%%%%%%%%%%%%%%%%%%%%%%%%%%
%%%%%%%%%%%%%%%%%%%%%%%%%%%%%%%%%%%%%%%%%%%%%%%%%%%%%%%%%%%%%%%%%%%%%%%%%%%%%%%%%%%%%%%%%%%%%%%%%%%%%%%%%%%%%%%%%%%%%%%%%
\section{Proofs of Lemma \ref{Lemma_5.1} and Lemma \ref{Lemma_5.2}}
In this section we continue use notation $f = O(g)$ meaning that $|f|\leq c|g|$ for some 
$c>0$ depending only on $r_1, r_2, c_1,c_2, a, b$ and $\varepsilon$.

Let
\begin{equation} 
\vec{\theta}^{(k)} = (\theta_1,\ldots,\theta_{k-1},0,\theta_{k+1}, \ldots, \theta_n).
\end{equation}
\begin{Proof} { \it Lemma \ref{Lemma_5.1}.}
Using (\ref{3.9}), (\ref{5.1}) and (\ref{5.4}), we find that
\begin{equation}\label{6.1}
	|<T>_{\mathbb{R}^n \setminus \Omega}| 
	\leq
	\int\limits_{\mathbb{R}^n \setminus \Omega}
	P(\|\vec{\theta}\|_\infty^s) e^{-\vec{\theta}^T A \vec{\theta} } d\vec{\theta} = O\left(\exp(-c_6n^{2\varepsilon})\right)<1>.
\end{equation}
For simplicity, let $k =1$. Using (\ref{6.1}), we get that
\begin{equation}
	\begin{aligned}
		<\phi_1^p(\vec{\theta})\tilde{T}(\vec{\theta})>_{F, \Omega} =
		\int\limits_{\Omega} \phi_1^p(\vec{\theta})\tilde{T}(\theta_2,\ldots,\theta_n)
		 e^{-\vec{\theta}^T A \vec{\theta} + H(\vec{\theta})} d\vec{\theta}
		=\\ = 	\int\limits_{ U_n(r_2n^{\varepsilon})} \phi_1^p(\vec{\theta})\tilde{T}(\theta_2,\ldots,\theta_n)
		 e^{-\vec{\theta}^T A \vec{\theta} + H(\vec{\theta})} d\vec{\theta} + O\left(\exp(-c_6n^{2\varepsilon})\right)<1>,\\
		 p=1,2,3,4.
	\end{aligned}
\end{equation}
%where we let $H(\vec{\theta}) \equiv 0$ for $\vec{\theta} \in \mathbb{R}^n \setminus \Omega$.
Combining (\ref{R_2}) and the Mean Value Theorem, we find that  for $\vec{\theta}\in U_n(r_2n^{\varepsilon})$
\begin{equation}\label{6.4}
	H(\vec{\theta}) - H(\vec{\theta}^{(1)}) = O(n^{-1+4\varepsilon}).
\end{equation}
Using (\ref{5.2}), we get that
\begin{equation}\label{6.5}
	\begin{aligned}
			\int\limits_{ U_n(r_2n^{\varepsilon})} \phi_1^p(\vec{\theta})\tilde{T}(\theta_2,\ldots,\theta_n)
		 e^{-\vec{\theta}^T A \vec{\theta} + H(\vec{\theta})} d\vec{\theta} + O\left(\exp(-c_6n^{2\varepsilon})\right)<1> = \\
		 = \int\limits_{-r_2n^{\varepsilon}}\limits^{r_2n^{\varepsilon}}\cdots \int\limits_{-r_2n^{\varepsilon}}\limits^{r_2n^{\varepsilon}}
	 \tilde{T}e^{-\,g_1(\theta_2,\ldots,\theta_{n}) + H(\vec{\theta}^{(1)})} \\
	\left( \int \limits_{-r_2n^{\varepsilon}}\limits^{r_2n^{\varepsilon}} \phi_1^p e^{-\phi_1^2(\vec{\theta}) + H(\vec{\theta}) - H(\vec{\theta}^{(1)})} d \theta_1  \right)
	d \theta_2 \ldots d \theta_n,\\ p =0, 1,2,3,4.
	\end{aligned}
\end{equation}
Combining (\ref{3.9}) and (\ref{6.4}), we find that for $\vec{\theta}^{(1)}\in U_n(r_2n^{\varepsilon})$
\begin{equation}\label{6.6}
	\begin{aligned}
	\int \limits_{-r_2n^{\varepsilon}}\limits^{r_2n^{\varepsilon}} &\phi_1^p e^{-\phi_1^2(\vec{\theta}) + H(\vec{\theta}) - 		
	H(\vec{\theta}^{(1)})} d \theta_1 = \\= 
	 \int\limits_{-\infty}\limits^{+\infty} \phi_1^p &e^{-\phi_1^2(\vec{\theta})}  d \theta_1 
 + 	O\left(\exp(-c_7n^{2\varepsilon})\right)\int\limits_{-\infty}\limits^{+\infty}  e^{-\phi_1^2(\vec{\theta})}  d \theta_1 +\\
 &+ \int \limits_{|\phi_1(\vec{\theta})|\leq r_3n^{\varepsilon}} \phi_1^p e^{-\phi_1^2(\vec{\theta})} \left(e^{H(\vec{\theta}) - 		
	H(\vec{\theta}^{(1)})}-1\right) d \theta_1, \  p =0,1,2,3,4,
	\end{aligned}
\end{equation}
where $c_7= c_7(r_2,c_1,c_2,a,b,\varepsilon)>0$, $r_3= r_3(r_2,c_1,c_2,a,b,\varepsilon)>0.$ 

For $p = 2,4$, we have that:
\begin{equation}\label{6.7}
	\begin{aligned}
	\int\limits_{-\infty}\limits^{+\infty} \phi_1^2 e^{-\phi_1^2(\vec{\theta})}  d \theta_1 = \frac{1}{2}
	\int\limits_{-\infty}\limits^{+\infty}  e^{-\phi_1^2(\vec{\theta})}  d \theta_1,\\
	\int\limits_{-\infty}\limits^{+\infty} \phi_1^4 e^{-\phi_1^2(\vec{\theta})}  d \theta_1 = \frac{3}{4}
	\int\limits_{-\infty}\limits^{+\infty}  e^{-\phi_1^2(\vec{\theta})}  d \theta_1,
	\end{aligned}
\end{equation}
\begin{equation}\label{6.8}
	\begin{aligned}
	\int \limits_{|\phi_1(\vec{\theta})|<r_3n^{\varepsilon}} \phi_1^p e^{-\phi_1^2(\vec{\theta})} \left(e^{H(\vec{\theta}) - 		
	H(\vec{\theta}^{(1)})}-1\right) d \theta_1 = \\ =
	O(n^{-1+4\epsilon})\int\limits_{-\infty}\limits^{+\infty} \phi_1^p  e^{-\phi_1^2(\vec{\theta})}  d \theta_1 = \\
	=O(n^{-1+4\epsilon})\int\limits_{-\infty}\limits^{+\infty}  e^{-\phi_1^2(\vec{\theta})}  d \theta_1,
	\\ \text{for } \vec{\theta}^{(1)}\in U_n(r_2n^{\varepsilon}), \ \ \  p = 0, 2,4.
	\end{aligned}
\end{equation}
Combining (\ref{6.1})-(\ref{6.8}), we obtain (\ref{phi^2}) and (\ref{phi^4}).

For $p=1,3$, we have that:
\begin{equation}\label{6.9}
	\int\limits_{-\infty}\limits^{+\infty} \phi_1^p e^{-\phi_1^2(\vec{\theta})}  d \theta_1 = 0
\end{equation}
\begin{equation}\label{6.10}
	\begin{aligned}
	\int \limits_{|\phi_1(\vec{\theta})|\leq r_3n^{\varepsilon}} \phi_1^p &e^{-\phi_1^2(\vec{\theta})} \left(e^{H(\vec{\theta}) - 		
	H(\vec{\theta}^{(1)})}-1\right) d \theta_1 =\\ &= \int \limits_{0 \leq \phi_1(\vec{\theta}) \leq r_3n^{\varepsilon}} |\phi_1|^p e^{-\phi_1^2(\vec{\theta})} \left(e^{H(\vec{\theta}) - 		
	H(\vec{\theta}^{(1)})}-1\right) d \theta_1 - \\ &- \int \limits_{r_3n^{\varepsilon} \leq \phi_1(\vec{\theta}) \leq 0} |\phi_1|^p e^{-\phi_1^2(\vec{\theta})} \left(e^{H(\vec{\theta}) - 		
	H(\vec{\theta}^{(1)})}-1\right) d \theta_1 =\\ &= 
	O(n^{-1+4\epsilon})\int\limits_{\phi_1\geq 0} |\phi_1|^p  e^{-\phi_1^2(\vec{\theta})}  d \theta_1 =\\
	&= O(n^{-1+4\epsilon})\int\limits_{-\infty}\limits^{+\infty}  e^{-\phi_1^2(\vec{\theta})}  d \theta_1, \\  
	&\ \ \ \ \text{ } \ \ \ \text{for } \vec{\theta}^{(1)}\in U_n(r_2n^{\varepsilon}), \ \ \ p = 1,3.
	\end{aligned}
\end{equation}
Combining (\ref{6.1})-(\ref{6.6}), \eqref{6.8}, (\ref{6.9}), (\ref{6.10}), we obtain (\ref{phi^2}) and (\ref{phi^4}).
\end{Proof}
%%%%%%%%%%%%%%%%%%%%%%%%%%%%%%%%%%%%%%%%%%%%%%%%%%%%%%%%%%%%%%%%%%%%%%%%%%%%%%%%%%%%%%%%%%%%%%%%%%%%%%%%%%%%%%%%%%%%%%%%%
%%%%%%%%%%%%%%%%%%%%%%%%%%%%%%%%%%%%%%%%%%%%%%%%%%%%%%%%%%%%%%%%%%%%%%%%%%%%%%%%%%%%%%%%%%%%%%%%%%%%%%%%%%%%%%%%%%%%%%%%%

\begin{Proof} { \it Lemma \ref{Lemma_5.2}.}
Let $T(\vec{\theta})$ satisfy 
\begin{equation}\label{5.12}
	\begin{aligned}
	|T(\vec{\theta})| = O(\|\vec{\theta}\|_\infty^s), \\
	\frac{\partial T(\vec{\theta})}{\partial \theta_k} = O(sn^{-1-\varepsilon})\sup_{\vec{\theta} \in \Omega}| T(\vec{\theta})|, 
	\ \ \vec{\theta} \in \Omega.
	\end{aligned}
\end{equation}
Combining (\ref{phi}) and (\ref{5.12}), we get that
\begin{equation}\label{5.13}
	\begin{aligned}
	<\phi_k(\vec{\theta})&T(\vec{\theta})>_{F, \Omega} = \\&=
	<\phi_k(\vec{\theta})T(\vec{\theta}^{(k)})>_{F, \Omega} + <\phi_k(\vec{\theta})(T(\vec{\theta})-T(\vec{\theta}^{(k)})>_{F, \Omega} 
	=\\ 
	&= <\phi_k(\vec{\theta})T(\vec{\theta}^{(k)})>_{F, \Omega}
	+ O(sn^{-1-\epsilon})<\sup_{\vec{\theta} \in \Omega}|\phi_k\theta_k T(\vec{\theta})|>_{F, \Omega} =\\
	&=O(sn^{-1+4\epsilon})<\sup_{\vec{\theta} \in \Omega}|T(\vec{\theta})|>_{F, \Omega} + O\left(\exp(-c_6n^{2\varepsilon})\right)<1>.
	\end{aligned}
\end{equation}
Using (\ref{A_ass}), we find that for $\vec{\theta}\in \Omega$
\begin{equation}\label{5.16}
%	\frac{\partial\left( M(\vec{\phi}(\vec{\theta})) - M(\vec{\theta})\right)}{\partial \theta_k} = 
	\frac{\partial M(\vec{\phi}(\vec{\theta}))}{\partial \theta_k} = O(sn^{(s-1)\varepsilon}) \sum\limits_{j : s_j \neq 0} 
	\frac{\partial \vec{\phi}_j (\vec{\theta})}{\partial \theta_k} = O(sn^{-1-\varepsilon}) \sup_{\vec{\theta} \in \Omega}
	|M(\vec{\phi}(\vec{\theta}))| .
\end{equation}
Combining (\ref{5.5}) and(\ref{5.13}), we  obtain (\ref{5.17})
\end{Proof}
%%%%%%%%%%%%%%%%%%%%%%%%%%%%%%%%%%%%%%%%%%%%%%%%%%%%%%%%%%%%%%%%%%%%%%%%%%%%%%%%%%%%%%%%%%%%%%%%%%%%%%%%%%%%%%%%%%%%%%%%%
%%%%%%%%%%%%%%%%%%%%%%%%%%%%%%%%%%%%%%%%%%%%%%%%%%%%%%%%%%%%%%%%%%%%%%%%%%%%%%%%%%%%%%%%%%%%%%%%%%%%%%%%%%%%%%%%%%%%%%%%%
%%%%%%%%%%%%%%%%%%%%%%%%%%%%%%%%%%%%%%%%%%%%%%%%%%%%%%%%%%%%%%%%%%%%%%%%%%%%%%%%%%%%%%%%%%%%%%%%%%%%%%%%%%%%%%%%%%%%%%%%%
%%%%%%%%%%%%%%%%%%%%%%%%%%%%%%%%%%%%%%%%%%%%%%%%%%%%%%%%%%%%%%%%%%%%%%%%%%%%%%%%%%%%%%%%%%%%%%%%%%%%%%%%%%%%%%%%%%%%%%%%%
%%%%%%%%%%%%%%%%%%%%%%%%%%%%%%%%%%%%%%%%%%%%%%%%%%%%%%%%%%%%%%%%%%%%%%%%%%%%%%%%%%%%%%%%%%%%%%%%%%%%%%%%%%%%%%%%%%%%%%%%%
%%%%%%%%%%%%%%%%%%%%%%%%%%%%%%%%%%%%%%%%%%%%%%%%%%%%%%%%%%%%%%%%%%%%%%%%%%%%%%%%%%%%%%%%%%%%%%%%%%%%%%%%%%%%%%%%%%%%%%%%%
%%%%%%%%%%%%%%%%%%%%%%%%%%%%%%%%%%%%%%%%%%%%%%%%%%%%%%%%%%%%%%%%%%%%%%%%%%%%%%%%%%%%%%%%%%%%%%%%%%%%%%%%%%%%%%%%%%%%%%%%%
%%%%%%%%%%%%%%%%%%%%%%%%%%%%%%%%%%%%%%%%%%%%%%%%%%%%%%%%%%%%%%%%%%%%%%%%%%%%%%%%%%%%%%%%%%%%%%%%%%%%%%%%%%%%%%%%%%%%%%%%%

\section{Proof of Proposition \ref{Proposition_3.1}}
In this section we use notation $f = O(g)$ meaning that $|f|\leq c|g|$ for some 
constant $c>0$ depending only on $\gamma$ and $\varepsilon$. 

The following lemma will be applied to estimate the determinant of a matrix close to
the identity matrix $I$.
\begin{Lemma}\label{Lemma_matrix1}
	Let $\|\|$ denote any matrix norm. Let $X$ be such an $n\times n$ matrix that $\left\|X\right\| < 1$. Then for fixed $m \geq 2$ 
	\begin{equation}
		\det(I+X) = \exp \left(  \sum\limits_{r=1}^{m-1} \frac{(-1)^{r+1}}{r}\, {\rm tr} (X^r) + E_m(X)   \right),
	\end{equation}
 	where ${\rm tr}(\cdot)$ is the trace function and
	\begin{equation}
		|E_m(X)|\leq \frac{n}{m}\,\frac{\left\|X\right\|^m}{1-\left\|X\right\|}.
	\end{equation}
\end{Lemma}
The proof of Lemma \ref{Lemma_matrix1} is based on estimation the trace of the matrix $\ln(I+X)$, using the representation as a convergent series. 
Lemma \ref{Lemma_matrix1} was also formulated and proved in\cite{Brendan1995}.

We have that
\begin{equation}\label{S_0000}
	S_0 = \frac{1}{n} \sum \limits_{r=1}\limits^n \int\limits_{V_0} 
 % \left( 
     \prod\limits_{\{v_j,v_k\}\in EG} \cos \Delta_{jk} \sum \limits_{T \in {\cal T}_r} 
     \prod \limits_{(v_j,v_k)\in ET} (1+i\tan \Delta_{jk})  
 % \right) 
 \
 d\vec{\xi},
\end{equation} 
where $\Delta_{jk} = \xi_j-\xi_k$ and
\begin{equation}
	\begin{aligned}
	V_0 = \{\vec{\xi}\in U_n(\pi/2) :\ |\xi_j - \xi_k|_\pi \leq n^{-1/2+\varepsilon} \text { for any} 
	1\leq j,k\leq n \},\\
	|\xi_j - \xi_k|_\pi = \min\limits_{l\in \mathbb{Z}} |\xi_j - \xi_k+ \pi l|.
	\end{aligned}
\end{equation}
Since the integrand is invariant under uniform translation of all
the $\xi_j$'s mod $\pi$, we can fix $\sum\limits_{k=1}\limits^n\xi_k = 0$  and multiply it by the ratio 
of its range $\pi$ to the length $n^{-1/2}$ of the vector $\frac{1}{n}[1,1,\ldots,1]^T$. 
Thus, for any $1\leq r\leq n$, we get that
\begin{equation}\label{8.5}
\begin{aligned}
S_0 = \pi n^{1/2}  \int\limits_{L \cap V_0\cap U_n(n^{-1/2+\varepsilon})}
%\left( 
   T(\vec{\xi})
% \right) 
  \ dL,\\
  T(\vec{\xi}) = \frac{1}{n} \sum \limits_{r=1}\limits^n    
  \prod\limits_{(v_j,v_k)\in EG} \cos \Delta_{jk} \sum \limits_{T \in {\cal T}_r} 
     \prod \limits_{(v_j,v_k)\in ET} (1+i\tan \Delta_{jk}), 
 \end{aligned}
\end{equation}
where $L$ denotes the orthogonal complement to the vector $[1,1,\ldots,1]^T$. 

Let define $n\times n$ matrix $B$ by
	\begin{equation}
		B_{jk} = 
		\left\{
		\begin{array}{ll}
			-\tan\Delta_{jk},& \text{ for } \{v_j, v_k\} \in EG,\\
			\sum\limits_{l: (v_j, v_l) \in EG} \tan \Delta_{jl},& \text{ for } k=j,\\
			0 & \text{ otherwise }.
			\end{array}
		\right.
	\end{equation}
	Using Theorem \ref{Tutte} with the matrix $Q + iB$, we get that 
	\begin{equation}\label{Lemma_using_Tutte_000}
			\sum\limits_{r=1}\limits^{n}\sum\limits_{T \in {\cal T}_r} \prod\limits_{(v_j,v_k)\in ET} (1+i\tan \Delta_{jk}) = 
			\sum\limits_{r=1}\limits^{n} M_r,	
	\end{equation} 	
	where $M_r$ denotes the
	principal minor of $A$ formed by removing row $r$ and column $r$.
	Since the vector $[1,1,\ldots,1]^T$ is the common eigenvector of the matrices $Q$ and $B$, corresponding to the eigenvalue $0$, we find that
	\begin{equation}\label{Lemma_using_Tutte_001}
		\sum\limits_{r=1}\limits^{n} M_r = \frac{\det(\hat{Q} + iB)}{n},
	\end{equation}
	where $\hat{Q} = Q + J$ and  $J$ denotes the matrix with every entry $1$. 
	Note that 
\begin{equation}\label{Lemma_using_Tutte_Delta}
		|\Delta_{jk}|  \leq  n^{-1/2+\varepsilon}, \ \ \ \vec{\xi} \in V_0\cap U_n(n^{-1/2+\varepsilon}),
	\end{equation}
			\begin{equation}\label{Lemma_using_Tutte_B}
		||B||_1 = \max_{j}{\sum\limits_{k=1}^{n}} |B_{jk}| = O(n^{1/2+\varepsilon}),\ \ \ \vec{\xi} \in U_n(n^{-1/2+\varepsilon}).
	\end{equation}
	Let $\Phi = B\hat{Q}^{-1}$. Using \eqref{3.9}, \eqref{4.3}, \eqref{5.29} and \eqref{Lemma_using_Tutte_B}, we get that
	\begin{equation}
		||\Phi||_1 \leq ||B||_1 ||\hat{Q}^{-1}||_1  = O(n^{-1/2+\varepsilon}),\ \ \ \vec{\xi} \in U_n(n^{-1/2+\varepsilon}).
	\end{equation}
	 Using Lemma \ref{Lemma_matrix1} with the matrix $i\Phi$, we find that 
	  \begin{equation}\label{Lemma_using_Tutte_I+Phi}
	 	\det(I+i \Phi) = \exp \left( i \mbox{tr}  \Phi + \frac{\mbox{tr} \Phi^2}{2}  + O(n^{-1/2+3\varepsilon}) \right),\ \ \ \vec{\xi} \in U_n(n^{-1/2+\varepsilon}).
	 \end{equation}
	 Let 
	 \begin{equation}
	 	B = B_{skew} + B_{diag},
	 \end{equation}
	 where $B_{skew}$ is the skew-symmetric matrix and $B_{diag}$ is the diagonal matrix.
	 Since $\hat{Q}$ is the symmetric matrix 
 	 \begin{equation}\label{Bskewsmall}
 	  	\mbox{tr}(B_{skew}\hat{Q}^{-1}) = 0.
 	 \end{equation}
 	Using (\ref{Lemma_using_Tutte_Delta}), we find that
 	\begin{equation}\label{Lemma_using_Tutte_DeltaB-L}
 	 		||B_{diag} - \Lambda||_2  = O(n^{-1/2+3\varepsilon}),\ \ \ \vec{\xi} \in U_n(n^{-1/2+\varepsilon}),
 	 \end{equation}
 	 where $\Lambda$ denotes the diagonal matrix whose diagonal elements are equal to the components of the vector $Q\vec{\theta}$.
 	 Combining \eqref{4.2} and \eqref{Lemma_using_Tutte_DeltaB-L}, we get that 
 	 \begin{equation}\label{Bdiag-lambda}
 	 	\begin{aligned}
 	 		\left|\mbox{tr}\left((B_{diag} - \Lambda)\hat{Q}^{-1}\right)\right| \leq n ||B_{diag} - \Lambda||_2||\hat{Q}^{-1}||_2 = O(n^{-1/2+3\varepsilon}),
 	 		 \\ \vec{\xi} \in U_n(n^{-1/2+\varepsilon}).
 	 	\end{aligned}
 	 \end{equation}
 	 Using (\ref{Bskewsmall}) and (\ref{Bdiag-lambda}), we obtain that 
 	 \begin{equation}\label{Lemma_using_Tutte_Phi}
 	 	\begin{aligned}
 	 		\mbox{tr}\Phi = \mbox{tr}(B_{diag}\hat{Q}^{-1}) = \mbox{tr}(\Lambda\hat{Q}^{-1})+ O(n^{-1/2+3\varepsilon}) = \\ = 
 	 		\vec{\xi}^T Q \vec{\alpha} + O(n^{-1/2+3\varepsilon}), \ \ \ \vec{\xi} \in V_0,
 	 		\end{aligned}
 	 \end{equation} 
 	 where $\vec{\alpha}$ is the vector composed of the diagonal elements of the matrix $\hat{Q}^{-1}$.
 	 
 	 Using the property of the trace function 
 	 \begin{equation}\label{XY=YX}
 	 		\mbox{tr}(XY) = \mbox{tr}(YX),
 	 \end{equation} we have that
 	 \begin{equation}\label{Lemma_using_Tutte_Phi2}
 	 		\mbox{tr} \Phi^2 = \mbox{tr}(B_{skew}\hat{Q}^{-1})^2 + \mbox{tr}(B_{diag}\hat{Q}^{-1})^2 + 
						2\, \mbox{tr}\left(B_{skew}\hat{Q}^{-1}B_{diag}\hat{Q}^{-1}\right).
 	 \end{equation}
 	 Since  $B_{skew}$ is the skew-symmetric matrix and $\hat{Q}^{-1}B_{diag}\hat{Q}^{-1}$ is the symmetric matrix, we find that
 	 \begin{equation}\label{Lemma_using_Tutte_B0}
 	 		\mbox{tr}\left(B_{skew}\hat{Q}^{-1}B_{diag}\hat{Q}^{-1}\right) = 0.
 	 \end{equation}
 	 According to \eqref{eq5.17}, we have that
 	 \begin{equation} \label{HS}
 	 	\begin{aligned}
 	 		\mbox{tr} X^2 &\leq ||X||_{HS}^2, \\
 	 		||XY||_{HS} &\leq ||X||_{HS}||Y^T||_2.
 	 	\end{aligned}
 	 \end{equation}	
 	 	Therefore we get that
 	 \begin{equation}\label{Lemma_using_Tutte_trBsk}
 	 		\left|\mbox{tr}(B_{skew}\hat{Q}^{-1})^2\right| \leq ||B_{skew}\hat{Q}^{-1}||^2_{HS}. 
 	 \end{equation}
 	 Combining \eqref{4.2} and (\ref{Lemma_using_Tutte_Delta}), we obtain that 
 	 \begin{equation}\label{Lemma_using_Tutte_Bsk}
 	 	||B_{skew}\hat{Q}^{-1}||_{HS} \leq ||\hat{Q}^{-1}||_2||B_{skew}||_{HS} 
 	 	 = O(n^{-1/2 +\varepsilon}), \ \ \ \vec{\xi} \in V_0.
 	 \end{equation}
 	 Using  \eqref{4.2}, (\ref{Lemma_using_Tutte_Delta}) and (\ref{Lemma_using_Tutte_DeltaB-L}), we get that
 	 \begin{equation}\label{8.24}
 	 	\begin{aligned}
 	 		\left|\mbox{tr}\left( (B_{diag}-\Lambda)\hat{Q}^{-1}B_{diag}\hat{Q}^{-1}\right)\right| &\leq \\ \leq 
 	 		n \|\hat{Q}^{-1}\|_2^2 ||(B_{diag}-\Lambda)||_2&||B_{diag}||_2 = O(n^{-1+4\varepsilon}), \ \ \ \vec{\xi} \in V_0,
 	 	\end{aligned}
 	 \end{equation}
 	 and
 	 \begin{equation}\label{8.25}
 	 	\begin{aligned}
 	 		\left|\mbox{tr}\left( (B_{diag}-\Lambda)\hat{Q}^{-1}(B_{diag}-\Lambda)\hat{Q}^{-1}\right)\right| &\leq \\ \leq
 	 		n \|\hat{Q}^{-1}\|_2^2 ||(B_{diag}-\Lambda)||^2_2 &= O(n^{-2+6\varepsilon}), \ \ \ \vec{\xi} \in V_0.
 	 	\end{aligned}
 	 \end{equation}
 	 Combining \eqref{XY=YX}, \eqref{8.24} and \eqref{8.25}, we obtain that
 	 \begin{equation}\label{Lemma_using_Tutte_Bdiag2}
 	 		 \mbox{tr}(B_{diag}\hat{Q}^{-1})^2 =  \mbox{tr}(\Lambda\hat{Q}^{-1})^2+ O(n^{-1+4\varepsilon}), \ \ \ \vec{\xi} \in V_0.
 	 \end{equation}
 	 Combining (\ref{Lemma_using_Tutte_Phi2}), (\ref{Lemma_using_Tutte_B0}), (\ref{Lemma_using_Tutte_trBsk}), (\ref{Lemma_using_Tutte_Bsk}) 
 	 and (\ref{Lemma_using_Tutte_Bdiag2}), we obtain that
 	 \begin{equation}\label{Lemma_using_Tutte_Phi2f}
 	 		\mbox{tr} \Phi^2 = \mbox{tr}(\Lambda\hat{Q}^{-1})^2 + O(n^{-1 +4\varepsilon}),\ \ \ \vec{\xi} \in V_0.
 	 \end{equation}
 	   Using (\ref{Lemma_using_Tutte_Phi}) and (\ref{Lemma_using_Tutte_Phi2f}) in (\ref{Lemma_using_Tutte_I+Phi}), we get that
 	  \begin{equation}\label{Lemma_using_Tutte_003}
 	  	\det(I+i \Phi) = \exp \left( i\,\vec{\theta}^T Q \vec{\alpha} + \frac{\mbox{tr}(\Lambda\hat{Q}^{-1})^2}{2}+
 	  	 O(n^{-1/2+4\varepsilon}) \right),\ \ \ \vec{\xi} \in V_0.
 	  \end{equation}
 	 
 	 	By Taylor's theorem  we have that for $\vec{\xi} \in V_0$
\begin{equation}\label{Lemma_SO_Taylor}
		\begin{aligned}
		\prod\limits_{(v_j,v_k)\in EG} &\cos \Delta_{jk} = \\ =
		\exp&\left(-\frac{1}{2}\sum\limits_{(v_j,v_k)\in EG}\Delta_{jk}^2 - 
			\frac{1}{12}\sum\limits_{(v_j,v_k)\in EG}\Delta_{jk}^4 + O(n^{-1+6\varepsilon})\right).
		\end{aligned}
	\end{equation}
		Note also that
	\begin{equation}\label{8.30}
		\sum\limits_{(v_j,v_k)\in EG}\Delta_{jk}^2 = \vec{\xi}^T Q \vec{\xi}.
	\end{equation}
 	 
 	 Putting together \eqref{8.5} for $r=1, 2\ldots,n$, \eqref{Lemma_using_Tutte_000}, \eqref{Lemma_using_Tutte_001}, \eqref{Lemma_using_Tutte_003}, \eqref{Lemma_SO_Taylor}, \eqref{8.30},	 we obtain that
\begin{equation}\label{8.31}
		S_0  = \pi n^{-3/2} \det \hat{Q} \left( \mbox{Int}' +  O\left(n^{-1/2+6\epsilon}\right) \mbox{Int}''\right),
	\end{equation}
	where 
		\begin{equation}\label{8.32}
		\begin{aligned}
		\mbox{Int}' =\int\limits_{ L \cap V_0 \cap U_n(n^{-1/2+\varepsilon})}	
		\exp(
				i\,\vec{\xi}^T Q \vec{\alpha} + F(\vec{\xi})) 
		 dL,\\
		 \mbox{Int}'' =\int\limits_{ L \cap V_0 \cap U_n(n^{-1/2+\varepsilon})}	
		e^{				F(\vec{\xi})}
		 dL,\\
		 F(\xi) = - \frac{1}{2} \vec{\xi}^T \hat{Q} \vec{\xi} 	- \frac{1}{12}\sum\limits_{\{v_j,v_k\}\in EG} \Delta_{jk}^4 + \frac{1}{2}R(\vec{\xi}).
		 \end{aligned}
	\end{equation} 
where $R(\vec{\xi}) = \mbox{tr} (\Lambda(\vec{\xi}) \hat{Q}^{-1} \Lambda(\vec{\xi}) \hat{Q}^{-1})$.

	Let ${\rm Pr}(\vec{\xi)}$ be the orthogonal projection  $\vec{\xi}$ onto the space $L$, where
$L$ denotes the orthogonal complement to the vector $[1,1,\ldots,1]^T$. 
Note that 
	\begin{equation}
		{\rm Pr}(\vec{\xi})= \vec{\xi} 
		- \bar{\xi} [1,1,\ldots,1]^T,
	\end{equation}
	where 
	\begin{equation}
 		\bar{\xi} = \frac{\xi_1 + \xi_2 + \ldots  \xi_n}{n}.
	\end{equation}
	Thus
	\begin{equation}\label{U_n_in_P}
		U_n( \frac{1}{2}n^{-1/2+\varepsilon}) \subset \left\{ \vec{\xi} \ : \ {\rm Pr}(\vec{\xi}) \in V_0\cap U_n(n^{-1/2+\varepsilon})\right\}
	\end{equation}
We also note that
\begin{equation}\label{P_Q}
	Q\vec{\xi} = Q{\rm Pr}(\vec{\xi}).
\end{equation}
Therefore the integrand of \eqref{8.32} does not change under the substitution of vector $\vec{\xi}$ by 
vector ${\rm Pr}(\vec{\xi})$ and
\begin{equation}
	\begin{aligned}
		\mbox{Int}'
		 =  \int\limits_{{\rm Pr}(\vec{\xi}) \in  L \cap V_0\cap U_n(n^{-1/2+\varepsilon})}	
		e^{
				i\,\vec{\xi}^T Q \vec{\alpha} + F(\vec{\xi})} 
		 d\vec{\xi} \Bigg/
		 \int\limits_{-\infty}\limits^{+\infty} e^{-\frac{1}{2}nx^2} dx =\\= 
		 \frac{n^{1/2}}{\sqrt{2\pi}}
		 \int\limits_{{\rm Pr}(\vec{\xi}) \in V_0\cap U_n(n^{-1/2+\varepsilon})}	
		e^{
				i\,\vec{\xi}^T Q \vec{\alpha} + F(\vec{\xi})} 
		 d\vec{\xi} ,
	\end{aligned}
\end{equation}

\begin{equation}
	\begin{aligned}
		\mbox{Int}''
		 =  \int\limits_{{\rm Pr}(\vec{\xi}) \in  L \cap V_0\cap U_n(n^{-1/2+\varepsilon})}	
		e^{
				 F(\vec{\xi})} 
		 d\vec{\xi} \Bigg/
		 \int\limits_{-\infty}\limits^{+\infty} e^{-\frac{1}{2}nx^2} dx =\\= 
		 \frac{n^{1/2}}{\sqrt{2\pi}}
		 \int\limits_{{\rm Pr}(\vec{\xi}) \in V_0\cap U_n(n^{-1/2+\varepsilon})}	
		e^{
				F(\vec{\xi})} 
		 d\vec{\xi} ,
	\end{aligned}
\end{equation}

Using notations \eqref{eq5.6} - \eqref{eq5.11}, and formulas \eqref{theta^0},  \eqref{U_n_in_P}, \eqref{eq5.67}, we get that
\begin{equation}\label{8.39}
	\begin{aligned}
	\mbox{Int}' = \left(1 + O\left(\exp(-c_7n^{2\varepsilon})\right)\right)\frac{n^{1/2}}{\sqrt{2\pi}}\, \mbox{Int}, \ \ \ \
	\mbox{Int}'' = O(1)\frac{n^{1/2}}{\sqrt{2\pi}}\,\mbox{Int},
	\end{aligned}
\end{equation}
where $c_7>0$ depends only on $\gamma$ and $\varepsilon$.  

Combining \eqref{8.31} and \eqref{8.39}, we obtain \eqref{S_0+++}.

%%%%%%%%%%%%%%%%%%%%%%%%%%%%%%%%%%%%%%%%%%%%%%%%%%%%%%%%%%%%%%%%%%%%%%%%%%%%%%%%%%%%%%%%%%%%%%%%%%%%%%%%%%%%%%%%%%%%%%%%%
%%%%%%%%%%%%%%%%%%%%%%%%%%%%%%%%%%%%%%%%%%%%%%%%%%%%%%%%%%%%%%%%%%%%%%%%%%%%%%%%%%%%%%%%%%%%%%%%%%%%%%%%%%%%%%%%%%%%%%%%%
%%%%%%%%%%%%%%%%%%%%%%%%%%%%%%%%%%%%%%%%%%%%%%%%%%%%%%%%%%%%%%%%%%%%%%%%%%%%%%%%%%%%%%%%%%%%%%%%%%%%%%%%%%%%%%%%%%%%%%%%%
%%%%%%%%%%%%%%%%%%%%%%%%%%%%%%%%%%%%%%%%%%%%%%%%%%%%%%%%%%%%%%%%%%%%%%%%%%%%%%%%%%%%%%%%%%%%%%%%%%%%%%%%%%%%%%%%%%%%%%%%%
%%%%%%%%%%%%%%%%%%%%%%%%%%%%%%%%%%%%%%%%%%%%%%%%%%%%%%%%%%%%%%%%%%%%%%%%%%%%%%%%%%%%%%%%%%%%%%%%%%%%%%%%%%%%%%%%%%%%%%%%%
%%%%%%%%%%%%%%%%%%%%%%%%%%%%%%%%%%%%%%%%%%%%%%%%%%%%%%%%%%%%%%%%%%%%%%%%%%%%%%%%%%%%%%%%%%%%%%%%%%%%%%%%%%%%%%%%%%%%%%%%%
%%%%%%%%%%%%%%%%%%%%%%%%%%%%%%%%%%%%%%%%%%%%%%%%%%%%%%%%%%%%%%%%%%%%%%%%%%%%%%%%%%%%%%%%%%%%%%%%%%%%%%%%%%%%%%%%%%%%%%%%%
%%%%%%%%%%%%%%%%%%%%%%%%%%%%%%%%%%%%%%%%%%%%%%%%%%%%%%%%%%%%%%%%%%%%%%%%%%%%%%%%%%%%%%%%%%%%%%%%%%%%%%%%%%%%%%%%%%%%%%%%%
\section*{Acknowledgements}
 This work was carried out under the supervision of S.P. Tarasov and supported in part by
RFBR grant no 11-01-00398a.
%This work was supported by the RFBR grant .

\noindent
{ {\bf M.I. Isaev}\\
Moscow Institute of Physics and Technology,

141700 Dolgoprudny, Russia\\
Centre de Math\'ematiques Appliqu\'ees, Ecole Polytechnique,

91128 Palaiseau, France\\
e-mail: \tt{Isaev.M.I@gmail.com}}\\

\end{document}